\documentclass[12pt,a4paper]{amsart}

\usepackage{amsmath,amssymb,amscd}
\font\fiverm=cmr5 % Now needed to avoid an error when inputing pictex
\input{prepictex}
\input{pictex}
\input{postpictex}

\theoremstyle{plain}
\newtheorem{thm}{Theorem}[section]
\newtheorem{prop}[thm]{Proposition}
\newtheorem{lemma}[thm]{Lemma}
\newtheorem{cor}[thm]{Corollary}
\newtheorem{defn}[thm]{Definition}
 \newcommand{\A}{{\mathcal A}}
        \newcommand{\D}{{\mathcal D}}\newcommand{\HH}{{\mathcal H}}
        \newcommand{\LL}{{\mathcal L}}
        \newcommand{\B}{{\mathcal B}}
        \newcommand{\K}{{\mathcal K}}

        \newcommand{\s}{\sigma}
\newcommand{\al}{\alpha}

        \newcommand{\dd}{|\D|}
        \newcommand{\n}{\parallel}
\newcommand{\bma}{\left(\begin{array}{cc}}
\newcommand{\ema}{\end{array}\right)}
\newcommand{\bca}{\left(\begin{array}{c}}
\newcommand{\eca}{\end{array}\right)}
\def\clsp{\overline{\operatorname{span}}}

\def\Aut{\operatorname{Aut}}

\newcommand{\hideqed}{\renewcommand{\qed}{}}
        \newcommand{\R}{\mathbf R}
        \newcommand{\C}{\mathbf C}
        
\newcommand{\Z}{\mathbf Z}
\newcommand{\N}{\mathbf N}

\newcommand{\ben}{\begin{displaymath}}
        \newcommand{\een}{\end{displaymath}}
\newcommand{\be}{\begin{equation}}
\newcommand{\ee}{\end{equation}}

        \newcommand{\bean}{\begin{eqnarray*}}
        \newcommand{\eean}{\end{eqnarray*}}
\newcommand{\nno}{\nonumber\\}
\newcommand{\bea}{\begin{eqnarray}}
        \newcommand{\eea}{\end{eqnarray}}

\def\cross#1{\rlap{\hskip#1pt\hbox{$-$}}}
        
        \def\bigintcross{\cross{2.3}\int}

\newcommand{\nc}{\newcommand}
\nc{\nt}{\newtheorem} \nc{\gf}[2]{\genfrac{}{}{0pt}{}{#1}{#2}}
\nc{\mb}[1]{{\mbox{$ #1 $}}} \nc{\real}{{\mathbb R}}
\nc{\comp}{{\mathbb C}} \nc{\ints}{{\mathbb Z}}
\nc{\Ltoo}{\mb{L^2({\mathbf H})}} \nc{\rtoo}{\mb{{\mathbf R}^2}}
\nc{\slr}{{\mathbf {SL}}(2,\real)} \nc{\slz}{{\mathbf
{SL}}(2,\ints)} \nc{\su}{{\mathbf {SU}}(1,1)} \nc{\so}{{\mathbf
{SO}}} \nc{\hyp}{{\mathbb H}} \nc{\disc}{{\mathbf D}}
\nc{\torus}{{\mathbb T}}

\nc{\ca}{{\mathcal A}} \nc{\cag}{{{\mathcal A}^\Gamma}}
\nc{\cg}{{\mathcal G}} \nc{\chh}{{\mathcal H}} \nc{\ck}{{\mathcal
B}} \nc{\cl}{{\mathcal L}} \nc{\cm}{{\mathcal M}}
\nc{\cn}{{\mathcal N}} \nc{\cs}{{\mathcal S}} \nc{\cz}{{\mathcal
Z}}
\nc{\sind}{\sigma{\rm -ind}}
\newcommand{\la}{\langle}
\newcommand{\ra}{\rangle}

\begin{document}

\pagestyle{headings}

%\begin{frontmatter}



\title{The Noncommutative Geometry of Graph $C^*$-Algebras I: The Index
Theorem}\footnote{\noindent This research was supported by
Australian Research Council and a University of Newcastle Project
Grant}


\author{David Pask}
\email{david.pask@newcastle.edu.au}
\author{Adam Rennie}
\email{ adam.rennie@newcastle.edu.au}

%\corauth[cor1]{Corresponding author, Fax: 61 2 4921 6898}
\address{School of Mathematical and Physical Sciences\\
University of Newcastle, Callaghan\\ NSW Australia, 2308}

\begin{abstract}
We investigate conditions on a graph $C^*$-algebra for the
existence of a faithful semifinite trace. Using such a trace and
the natural gauge action of the circle on the graph algebra, we
construct a smooth $(1,\infty)$-summable semfinite spectral
triple. The local index theorem allows us to compute the pairing
with $K$-theory. This produces invariants in the $K$-theory of the
fixed point algebra, and these are invariants for a finer
structure than the isomorphism class of $C^*(E)$.

\vspace{3mm} \noindent {\bf Keywords:} Graph $C^*$-algebra, spectral
triple, index theorem, semifinite von Neumann algebra, trace,
$K$-theory, $KK$-theory.

\vspace{3mm} \noindent {\bf MSC (2000)} primary: 46L80, 58B34;
secondary 46L51, 46L08
\end{abstract}

%\begin{keyword}

%\end{keyword}
%\end{frontmatter}

\maketitle

% main text
\section{Introduction}
\label{intro} The aim of this paper, and the sequel \cite{PRen},
is to investigate the noncommutative geometry of graph
$C^*$-algebras. In particular we construct finitely summable
spectral triples to which we can apply the local index theorem.
The motivation for this is the need for new examples in
noncommutative geometry. Graph $C^*$-algebras allow us to treat a
large family of algebras in a uniform manner.

Graph $C^*$-algebras have been widely studied, see
\cite{BPRS,kpr,KPRR,H,PR,RSz,T} and the references therein. The
freedom to use both graphical and analytical tools make them
particularly tractable. In addition, there are many natural
generalisations of this family to which our methods will apply,
such as Cuntz-Krieger, Cuntz-Pimsner algebras, Exel-Laca algebras,
$k$-graph algebras and so on; for more information on these
classes of algebras see the above references and \cite{R}. We
expect these classes to yield similar examples.

One of the key features of this work is that the natural
construction of a spectral triple $(\A,\HH,\D)$ for a graph
$C^*$-algebra is almost never a spectral triple in the original
sense, \cite[Chapter VI]{C}. That is, the key requirement that for
all $a\in \A$ the operator $a(1+\D^2)^{-1/2}$ be a compact
operator on the Hilbert space $\HH$ is almost never true. However,
if we broaden our point of view to consider semifinite spectral
triples, where we require $a(1+\D^2)^{-1/2}$ to be in the ideal of
compact operators in a semifinite von Neumann algebra, we obtain
many $(1,\infty)$-summable examples. The only connected
$(1,\infty)$-summable example arising from our construction which
satisfies the original definition of spectral triples is the Dirac
triple for the circle.

{\bf The way we arrive at the correct notion of compactness is to
regard the fixed point subalgebra $F$ for the $S^1$ gauge action
on a graph algebra as the scalars.} This provides a unifying point
of view that will help the reader motivate the various
constructions, and understand the results. For instance the
$C^*$-bimodule we employ is a $C^*$-module over $F$, the range of
the ($C^*$-) index pairing lies in $K_0(F)$, the `differential'
operator $\D$ is linear over $F$ and it is the `size' of $F$ that
forces us to use a general semifinite trace. The single
$(1,\infty)$-summable example where the operator trace arises as
the natural trace is the circle, and in this case $F=\C$.

The algebras which arise from our construction, despite naturally
falling into the semifinite picture of spectral triples, are all
type I algebras, \cite{DHS}. Thus even when dealing with type I
algebras there is a natural and important role for general
semifinite traces.

Many of our examples arise from nonunital algebras. Fortunately,
graph $C^*$-algebras (and their smooth subalgebras) are
quasi-local in the sense of \cite{GGISV}, and many of the results
for smooth local algebras presented in \cite{R1,R2} are valid for
smooth quasi-local algebras. Here `local' refers to the
possibility of using a notion of `compact support' to deal with
analytical problems.

After some background material, we begin in Section \ref{triplesI}
by constructing an odd Kasparov module $(X,V)$ for $C^*(E)$-$F$,
where $F$ is the fixed point algebra. This part of the
construction applies to any locally finite directed graph with no
sources. The class $(X,V)$ can be paired with $K_1(C^*(E))$ to
obtain an index class in $K_0(F)$. This pairing is described in
the Appendix, and it is given  in terms of the index of Toeplitz
operators on the underlying $C^*$-module. We conjecture that this
pairing is the Kasparov product.

When our graph $C^*$-algebra has a faithful (semifinite, lower-semicontinuous)
 gauge invariant trace $\tau$, we can
define a canonical faithful (semifinite, lower semicontinuous)
trace $\tilde\tau$ on the endomorphism
algebra of the $C^*$-$F$-module $X$. Using $\tilde\tau$, in
Section \ref{triplesII} we construct a semifinite spectral triple
$(\A,\HH,\D)$ for a smooth subalgebra $\A\subset C^*(E)$.

The numerical index pairing of $(\A,\HH,\D)$ with $K_1(C^*(E))$
can be computed using the semifinite local index theorem,
\cite{CPRS2}, and we prove that
$$ \la K_1(C^*(E)),(\A,\HH,\D)\ra=\tilde\tau_*\la
K_1(C^*(E)),(X,V)\ra,$$ where $\la K_1(C^*(E)),(X,V)\ra\subset K_0(F)$
denotes the $K_0(F)$-valued index and $\tilde\tau_*$ is the map
induced on $K$-theory by $\tilde\tau$. We show by an example that
this pairing is an invariant of a finer structure than the
isomorphism class of $C^*(E)$.

To ensure that readers without a background in graph C*-algebras
or a background in spectral triples can access the results in this
paper, we have tried to make it self contained. The organisation
of the paper is as follows. Section \ref{background} describes
graph $C^*$-algebras and semifinite spectral triples, as well as
quasilocal algebras and the local index theorem. Section
\ref{traces} investigates which graph $C^*$-algebras have a
faithful positive trace, and we provide some necessary and some
sufficient conditions. In Section \ref{triplesI} we construct a
$C^*$-module for any locally finite graph $C^*$-algebra. Using the
generator of the gauge action on this $C^*$-module, we obtain a
Kasparov module whenever the graph has no sources, and so a
$KK$-class. In Section \ref{triplesII}, we restrict to those graph
$C^*$-algebras with a faithful gauge invariant trace, and
construct a spectral triple from our Kasparov module. Section
\ref{index} describes our results pertaining to the index theorem.

In the sequel to this paper, \cite{PRen}, we identify a large
subclass of our graph $C^*$-algebras with faithful trace which
satisfy a natural semifinite and nonunital generalisation of
Connes' axioms for noncommutative manifolds. These examples are
all one dimensional.

{\bf Acknowledgements} We would like to thank Iain Raeburn and
Alan Carey for many useful comments and support. We  also thank the referee
for many useful comments that have improved the work. In addition, we thank
Nigel
Higson for showing us a proof that the pairing in the Appendix does indeed
represent the Kasparov product.

\section{Graph $C^*$-Algebras and Semifinite Spectral Triples}\label{background}
\vspace{-7pt}
\subsection{The $C^*$-algebras of Graphs}\label{graphalg}
\vspace{-7pt}
 For a more detailed introduction to graph
$C^*$-algebras we refer the reader to \cite{BPRS,kpr} and the
references therein. A directed graph $E=(E^0,E^1,r,s)$ consists of
countable sets $E^0$ of vertices and $E^1$ of edges, and maps
$r,s:E^1\to E^0$ identifying the range and source of each edge.
{\bf We will always assume that the graph is} {\bf row-finite}
which means that each vertex emits at most finitely many edges.
Later we will also assume that the graph is \emph{locally finite}
which means  it is row-finite and each vertex receives at most
finitely many edges. We write $E^n$ for the set of paths
$\mu=\mu_1\mu_2\cdots\mu_n$ of length $|\mu|:=n$; that is,
sequences of edges $\mu_i$ such that $r(\mu_i)=s(\mu_{i+1})$ for
$1\leq i<n$.  The maps $r,s$ extend to $E^*:=\bigcup_{n\ge 0} E^n$
in an obvious way. A \emph{loop} in $E$ is a path $L \in E^*$ with
$s ( L ) = r ( L )$, we say that a loop $L$ has an exit if there
is $v = s ( L_i )$ for some $i$ which emits more than one edge. If
$V \subseteq E^0$ then we write $V \ge w$ if there is a path $\mu
\in E^*$ with $s ( \mu ) \in V$ and $r ( \mu ) = w$
(we also sometimes say
 that $w$ is downstream from $V$). A \emph{sink}
is a vertex $v \in E^0$ with $s^{-1} (v) = \emptyset$, a
\emph{source} is a vertex $w \in E^0$ with $r^{-1} (w) =
\emptyset$.

A \emph{Cuntz-Krieger $E$-family} in a $C^*$-algebra $B$ consists
of mutually orthogonal projections $\{p_v:v\in E^0\}$ and partial
isometries $\{S_e:e\in E^1\}$ satisfying the \emph{Cuntz-Krieger
relations}
\begin{equation*}
S_e^* S_e=p_{r(e)} \mbox{ for $e\in E^1$} \ \mbox{ and }\
p_v=\sum_{\{ e : s(e)=v\}} S_e S_e^*   \mbox{ whenever $v$ is not
a sink.}
\end{equation*}
%We shall typically use small letters $\{s_e,p_v\}$ for
%Cuntz-Krieger families in a $C^*$-algebra and large letters
%$\{S_e,P_v\}$ for Cuntz-Krieger families of operators on Hilbert
%space.

It is proved in \cite[Theorem 1.2]{kpr} that there is a universal
$C^*$-algebra $C^*(E)$ generated by a non-zero Cuntz-Krieger
$E$-family $\{S_e,p_v\}$.  A product
$S_\mu:=S_{\mu_1}S_{\mu_2}\dots S_{\mu_n}$ is non-zero precisely
when $\mu=\mu_1\mu_2\cdots\mu_n$ is a path in $E^n$. Since the
Cuntz-Krieger relations imply that the projections $S_eS_e^*$ are
also mutually orthogonal, we have $S_e^*S_f=0$ unless $e=f$, and
words in $\{S_e,S_f^*\}$ collapse to products of the form $S_\mu
S_\nu^*$ for $\mu,\nu\in E^*$ satisfying $r(\mu)=r(\nu)$ (cf.\
\cite[Lemma
  1.1]{kpr}).
Indeed, because the family $\{S_\mu S_\nu^*\}$ is closed under
multiplication and involution, we have
\begin{equation}
C^*(E)=\clsp\{S_\mu S_\nu^*:\mu,\nu\in E^*\mbox{ and
}r(\mu)=r(\nu)\}.\label{spanningset}
\end{equation}
The algebraic relations and the density of $\mbox{span}\{S_\mu
S_\nu^*\}$ in $C^*(E)$ play a critical role throughout the paper.
We adopt the conventions that vertices are paths of length 0, that
$S_v:=p_v$ for $v\in E^0$, and that all paths $\mu,\nu$ appearing
in (\ref{spanningset}) are non-empty; we recover $S_\mu$, for
example, by taking $\nu=r(\mu)$, so that $S_\mu S_\nu^*=S_\mu
p_{r(\mu)}=S_\mu$.

If $z\in S^1$, then the family $\{zS_e,p_v\}$ is another
Cuntz-Krieger $E$-family which generates $C^*(E)$, and the
universal property gives a homomorphism $\gamma_z:C^*(E)\to
C^*(E)$ such that $\gamma_z(S_e)=zS_e$ and $\gamma_z(p_v)=p_v$.
The homomorphism $\gamma_{\overline z}$ is an inverse for
$\gamma_z$, so $\gamma_z\in\Aut C^*(E)$, and a routine
$\epsilon/3$ argument using (\ref{spanningset}) shows that
$\gamma$ is a strongly continuous action of $S^1$ on $C^*(E)$. It
is called the \emph{gauge action}. Because $S^1$ is compact,
averaging over $\gamma$ with respect to normalised Haar measure
gives an expectation $\Phi$ of $C^*(E)$ onto the fixed-point
algebra $C^*(E)^\gamma$:
\[
\Phi(a):=\frac{1}{2\pi}\int_{S^1} \gamma_z(a)\,d\theta\ \mbox{ for
}\ a\in C^*(E),\ \ z=e^{i\theta}.
\]
The map $\Phi$ is positive, has norm $1$, and is faithful in the
sense that $\Phi(a^*a)=0$ implies $a=0$.

From Equation (\ref{spanningset}), it is easy to see that a graph
$C^*$-algebra is unital if and only if the underlying graph is
finite. When we consider infinite graphs, formulas which involve
sums of projections may contain infinite sums. To interpret these,
we use strict convergence in the multiplier algebra of $C^*(E)$:

\begin{lemma}\label{strict}
Let $E$ be a row-finite graph, let $A$ be a $C^*$-algebra
generated by a Cuntz-Krieger $E$-family $\{T_e,q_v\}$, and let
$\{p_n\}$ be a sequence of projections  in $A$. If $p_nT_\mu
T_\nu^*$ converges for every $\mu,\nu\in E^*$, then $\{p_n\}$
converges strictly to a projection $p\in M(A)$.
\end{lemma}

\begin{proof}
Since we can approximate any $a\in A=\pi_{T,q}(C^*(E))$ by a
linear combination of $T_\mu T_\nu^*$, an $\epsilon/3$-argument
shows that $\{p_na\}$ is Cauchy for every $a\in A$. We define
$p:A\to A$ by $p(a):=\lim_{n\to\infty}p_na$. Since
\[
b^*p(a)=\lim_{n\to\infty}b^*p_na=\lim_{n\to\infty}(p_nb)^*a=p(b)^*a,
\]
the map $p$ is an adjointable operator on the Hilbert $C^*$-module
$A_A$, and hence defines (left multiplication by) a multiplier $p$
of $A$ \cite[Theorem 2.47]{RW}. Taking adjoints shows that
$ap_n\to ap$ for all $a$, so $p_n\to p$ strictly. It is easy to
check that $p^2=p=p^*$.
\end{proof}

\vspace{-10pt}
\subsection{Semifinite Spectral Triples}
\vspace{-10pt}
 We begin with some semifinite versions of standard
definitions and results. Let $\tau$ be a fixed faithful, normal,
semifinite trace on the von Neumann algebra ${\mathcal N}$. Let
${\mathcal K}_{\mathcal N }$ be the $\tau$-compact operators in
${\mathcal N}$ (that is the norm closed ideal generated by the
projections $E\in\mathcal N$ with $\tau(E)<\infty$).

\begin{defn} A semifinite
spectral triple $(\A,\HH,\D)$ is given by a Hilbert space $\HH$, a
$*$-algebra $\A\subset \cn$ where $\cn$ is a semifinite von
Neumann algebra acting on $\HH$, and a densely defined unbounded
self-adjoint operator $\D$ affiliated to $\cn$ such that

1) $[\D,a]$ is densely defined and extends to a bounded operator
for all $a\in\A$

2) $a(\lambda-\D)^{-1}\in\K_\cn$ for all $\lambda\not\in{\R}\
\mbox{and all}\ a\in\A.$

3) The triple is said to be even if there is $\Gamma\in\cn$ such
that $\Gamma^*=\Gamma$, $\Gamma^2=1$,  $a\Gamma=\Gamma a$ for all
$a\in\A$ and $\D\Gamma+\Gamma\D=0$. Otherwise it is odd.
\end{defn}

\begin{defn}\label{qck} A semifinite spectral triple $(\A,\HH,\D)$ is $QC^k$ for $k\geq 1$
($Q$ for quantum) if for all $a\in\A$ the operators $a$ and
$[\D,a]$ are in the domain of $\delta^k$, where
$\delta(T)=[\dd,T]$ is the partial derivation on $\cn$ defined by
$\dd$. We say that $(\A,\HH,\D)$ is $QC^\infty$ if it is $QC^k$
for all $k\geq 1$.
\end{defn}

{\bf Note}. The notation is meant to be analogous to the classical
case, but we introduce the $Q$ so that there is no confusion
between quantum differentiability of $a\in\A$ and classical
differentiability of functions.

\noindent{\bf Remarks concerning derivations and commutators}.  By
partial derivation we mean that $\delta$ is defined on some
subalgebra of $\cn$ which need not be (weakly) dense in $\cn$.
More precisely, $\mbox{dom}\ \delta=\{T\in\cn:\delta(T)\mbox{ is
bounded}\}$. We also note that if $T\in{\mathcal N}$, one can show
that $[\dd,T]$ is bounded if and only if $[(1+\D^2)^{1/2},T]$ is
bounded, by using the functional calculus to show that
$\dd-(1+\D^2)^{1/2}$ extends to a bounded operator in $\cn$. In
fact, writing $\dd_1=(1+\D^2)^{1/2}$ and $\delta_1(T)=[\dd_1,T]$
we have \ben \mbox{dom}\ \delta^n=\mbox{dom}\ \delta_1^n\ \ \ \
\mbox{for all}\ n.\een We also observe that if $T\in\cn$ and
$[\D,T]$ is bounded, then $[\D,T]\in\cn$. Similar comments apply
to $[\dd,T]$, $[(1+\D^2)^{1/2},T]$. The proofs
can be found in \cite{CPRS2}.

The $QC^\infty$ condition places some restrictions on the algebras we
consider. Recall that a topological algebra is Fr\'{e}chet
if it is locally convex,
metrizable and complete, and that a subalgebra of a $C^*$-algebra is
a pre-$C^*$-algebra if it is stable under the holomorphic functional
calculus. For nonunital algebras, we consider only functions $f$ with $f(0)=0$.
\begin{defn}A $*$-algebra $\A$ is smooth if it is Fr\'{e}chet
and $*$-isomorphic to a proper dense subalgebra $i(\A)$ of a
$C^*$-algebra $A$ which is a pre-$C^*$-algebra.
\end{defn}  Asking for
$i(\A)$ to be a {\it proper} dense subalgebra of $A$ immediately
implies that the Fr\'{e}chet topology of $\A$ is finer than the
$C^*$-topology of $A$. We will denote the norm closure
$\overline{\A}=A$,  when  the norm closure $\overline{\A}$ is unambiguous.

If $\A$ is smooth in $A$ then $M_n(\A)$ is
smooth in $M_n(A)$, \cite{GVF,LBS}, so $K_*(\A)\cong K_*(A)$, the isomorphism
being induced by the inclusion map $i$. A
smooth algebra has a sensible spectral theory which agrees with that
defined using the $C^*$-closure, and the group of invertibles is
open.
The point of contact
between smooth algebras and $QC^\infty$ spectral triples is the
following Lemma, proved in \cite{R1}.

%We will always suppose that we can define the Fr\'{e}chet topology
%of $\A$ using a countable collection of submultiplicative
%seminorms which includes the $C^*$-norm of $\overline{\A}=A$, and
%note that the multiplication is jointly continuous, \cite{mal}. By
%replacing any seminorm $q$ by $\frac{1}{2}(q(a)+q(a^*))$, we may
%suppose that $q(a)=q(a^*)$ for all $a\in\A$.

%Stability under the holomorphic functional calculus extends to
%nonunital algebras, since the spectrum of an element in a
%nonunital algebra is defined to be the spectrum of this element in
%the  `one-point' unitization, though we must of course restrict to
%functions satisfying $f(0)=0$. Likewise, the definition of a
%Fr\'{e}chet algebra does not require a unit.

\begin{lemma}\label{smo} If $(\A,\HH,\D)$ is a $QC^\infty$ spectral triple, then
$(\A_\delta,\HH,\D)$ is also a $QC^\infty$ spectral triple, where
$\A_\delta$ is the completion of $\A$ in the locally convex
topology determined by the seminorms \ben
q_{n,i}(a)=\n\delta^nd^i(a)\n,\ \ n\geq 0,\ i=0,1,\een where
$d(a)=[\D,a]$. Moreover, $\A_\delta$ is a smooth algebra.
\end{lemma}

We call the topology on $\A$ determined by the seminorms $q_{n,i}$
of Lemma \ref{smo} the $\delta$-topology.

Whilst smoothness does not depend on whether $\A$ is unital or
not, many analytical problems arise because of the lack of a unit.
As in \cite{GGISV,R1,R2}, we make two definitions to address these
issues.

\begin{defn} An algebra $\A$ has local units if for every finite subset of
elements $\{a_i\}_{i=1}^n\subset\A$, there exists $\phi\in\A$ such
that for each $i$ \ben \phi a_i= a_i\phi=a_i.\een
\end{defn}

\begin{defn}
Let $\A$ be a Fr\'{e}chet algebra and $\A_c\subseteq\A$ be a dense
subalgebra with local units. Then we call  $\A$ a quasi-local
algebra (when $\A_c$ is understood.) If $\A_c$ is a dense ideal
with local units, we call $\A_c\subset\A$ local.
\end{defn}

Quasi-local algebras have an approximate unit $\{\phi_n\}_{n\geq
1}\subset\A_c$ such that for all $n$, $\phi_{n+1}\phi_n=\phi_n$,
\cite{R1}; we call this a local approximate unit.

{\bf Example} For a graph $C^*$-algebra $A=C^*(E)$, Equation
(\ref{spanningset}) shows that
$$ A_c=\mbox{span}\{S_\mu S_\nu^*:\mu,\nu\in E^*\ \mbox{and}\
r(\mu)=r(\nu)\}$$ is a dense subalgebra. It has local units
because
$$ p_{v}S_\mu S_\nu^*=\left\{\begin{array}{lr} S_\mu S_\nu^* &
v=s(\mu)\\ 0 & \mbox{otherwise}\end{array}\right..$$ Similar
comments apply to right multiplication by $p_{s(\nu)}$. By summing
the source and range projections (without repetitions) of all
$S_{\mu_i}S_{\nu_i}^*$ appearing in a finite sum
$$ a=\sum_ic_{\mu_i,\nu_i}S_{\mu_i}S_{\nu_i}^*$$
we obtain a local unit for $a\in A_c$. By repeating this process
for any finite collection of such $a\in A_c$ we see that $A_c$ has
local units.

We also require that when we have a spectral triple the operator
$\D$ is compatible with the quasi-local structure of the algebra,
in the following sense.

\begin{defn} If $(\A,\HH,\D)$ is a spectral triple, then we define $\Omega^*_\D(\A)$
to be the algebra generated by $\A$ and $[\D,\A]$.
\end{defn}

\begin{defn}\label{lst} A local spectral triple $(\A,\HH,\D)$ is a
spectral triple with $\A$ quasi-local such that there exists an
approximate unit $\{\phi_n\}\subset\A_c$ for $\A$ satisfying \ben
\Omega^*_\D(\A_c)=\bigcup_n\Omega^*_\D(\A)_n,\ \ {\rm where}\een
\ben
\Omega^*_\D(\A)_n=\{\omega\in\Omega^*_\D(\A):\phi_n\omega=\omega\phi_n=\omega\}.\een
\end{defn}

{\bf Remark} A local spectral triple has a local approximate unit
$\{\phi_n\}_{n\geq 1}\subset\A_c$  such that
$\phi_{n+1}\phi_n=\phi_n\phi_{n+1}=\phi_n$ and
$\phi_{n+1}[\D,\phi_n]=[\D,\phi_n]\phi_{n+1}=[\D,\phi_n]$, see \cite{R1,R2}.
We require this property to prove the summability
results we require.
\vspace{-5pt}
\subsection{Summability and the Local Index Theorem}
In the following, let $\mathcal N$ be a semifinite von Neumann
algebra with faithful normal trace $\tau$. Recall from \cite{FK}
that if $S\in\mathcal N$, the \emph{ t-th generalized singular
value} of S for each real $t>0$ is given by
$$\mu_t(S)=\inf\{||SE||\ : \ E \mbox{ is a projection in }
{\mathcal N} \mbox { with } \tau(1-E)\leq t\}.$$

The ideal $\LL^1({\mathcal N})$ consists of those operators $T\in
{\mathcal N}$ such that $\n T\n_1:=\tau( |T|)<\infty$ where
$|T|=\sqrt{T^*T}$. In the Type I setting this is the usual trace
class ideal. We will simply write $\LL^1$ for this ideal in order
to simplify the notation, and denote the norm on $\LL^1$ by
$\n\cdot\n_1$. An alternative definition in terms of singular
values is that $T\in\LL^1$ if $\|T\|_1:=\int_0^\infty \mu_t(T) dt
<\infty.$

Note that in the case where ${\mathcal N}\neq{\mathcal
B}({\mathcal H})$, $\LL^1$ is not complete in this norm but it is
complete in the norm $||.||_1 + ||.||_\infty$. (where
$||.||_\infty$ is the uniform norm). Another important ideal for
us is the domain of the Dixmier trace:
$${\mathcal L}^{(1,\infty)}({\mathcal N})=
\left\{T\in{\mathcal N}\ : \Vert T\Vert_{_{{\mathcal
L}^{(1,\infty)}}} :=   \sup_{t> 0}
\frac{1}{\log(1+t)}\int_0^t\mu_s(T)ds<\infty\right\}.$$

%There are related ideals for $p>1$: to describe them first set,
%for $p>1$,
%$$\psi_p(t)=\left\{\begin{array}{ll} t & \mbox{for } 0\leq t\leq 1\\
 %                                    t^{1-\frac{1}{p}} & \mbox{for } 1\leq t
%\end{array}\right..$$
%Then we define
%$${\mathcal L}^{(p,\infty)}({\mathcal N})
%=\left\{T\in{\mathcal N}\ : \Vert T\Vert_{_{{\mathcal
%L}^{(p,\infty)}}} := \sup_{t> 0}
%\frac{1}{\psi_p(t)}\int_0^t\mu_s(T)ds<\infty\right\}.$$ For $p>1$
%there is also the equivalent definition
%$${\mathcal L}^{(p,\infty)}({\mathcal N})
%=\left\{T\in{\mathcal N}\ : \sup_{t> 0}
%\frac{t}{\psi_p(t)}\mu_t(T)<\infty\right\}.$$ If
%$T\in\LL^{(p,\infty)}(\cn)$, then $T^p\in\LL^{(1,\infty)}(\cn)$.

We will suppress the $({\mathcal N})$ in our notation for these
ideals, as $\cn$ will always be clear from context. The reader
should note that ${\mathcal L}^{(1,\infty)}$ is often taken to
mean an ideal in the algebra $\widetilde{\mathcal N}$ of
$\tau$-measurable operators affiliated to ${\mathcal N}$, \cite{FK}. Our
notation is however consistent with that of \cite{C} in the
special case ${\mathcal N}={\mathcal B}({\mathcal H})$. With this
convention the ideal of $\tau$-compact operators, ${\mathcal
  K}({\mathcal N})$,
consists of those $T\in{\mathcal N}$ (as opposed to
$\widetilde{\mathcal N}$) such that \ben \mu_\infty(T):=\lim
_{t\to \infty}\mu_t(T)  = 0.\een
\begin{defn}\label{summable} A semifinite local spectral triple is
$(1,\infty)$-summable if \ben
a(\D-\lambda)^{-1}\in\LL^{(1,\infty)}\ \ \ \mbox{for all}\
a\in\A_c,\ \ \lambda\in\C\setminus\R.\een
Equivalently, $a(1+\D^2)^{-1/2}\in\LL^{(1,\infty)}$ for all $a\in \A_c$.
\end{defn}

{\bf Remark} If $\A$ is unital, $\ker\D$ is $\tau$-finite
dimensional. Note that the summability requirements are only for
$a\in\A_c$. We do not assume that elements of the algebra $\A$ are
all integrable in the nonunital case.

We need to briefly discuss the Dixmier trace, but fortunately we
will usually be applying it in reasonably simple situations. For
more information on semifinite Dixmier traces, see \cite{CPS2}.
For $T\in\LL^{(1,\infty)}$, $T\geq 0$, the function \ben
F_T:t\to\frac{1}{\log(1+t)}\int_0^t\mu_s(T)ds \een is bounded. For
certain generalised limits $\omega\in L^\infty(\R_*^+)^*$, we
obtain a positive functional on $\LL^{(1,\infty)}$ by setting
$$ \tau_\omega(T)=\omega(F_T).$$
 This is the
Dixmier trace associated to the semifinite normal trace $\tau$,
denoted $\tau_\omega$, and we extend it to all of
$\LL^{(1,\infty)}$ by linearity, where of course it is a trace.
The Dixmier trace $\tau_\omega$ is defined on the ideal
$\LL^{(1,\infty)}$, and vanishes on the ideal of trace class
operators. Whenever the function $F_T$ has a limit at infinity,
all Dixmier traces return the value of the limit. We denote the
common value of all Dixmier traces on measurable operators by
$\bigintcross$. So if $T\in\LL^{(1,\infty)}$ is measurable, for
any allowed functional $\omega\in L^\infty(\R_*^+)^*$ we have
$$\tau_\omega(T)=\omega(F_T)=\bigintcross T.$$

{\bf Example} Let $\D=\frac{1}{i}\frac{d}{d\theta}$ act on
$L^2(S^1)$. Then it is well known that the spectrum of $\D$
consists of eigenvalues $\{n\in\Z\}$, each with multiplicity one.
So, using the standard operator trace, the function
$F_{(1+\D^2)^{-1/2}}$ is
$$ N\to\frac{1}{\log 2N+1}\sum_{n=-N}^N(1+n^2)^{-1/2}$$
which is bounded. So $(1+\D^2)^{-1/2}\in\LL^{(1,\infty)}$ and for
any Dixmier trace $\mbox{Trace}_\omega$
$$\mbox{Trace}_\omega((1+\D^2)^{-1/2})=\bigintcross(1+\D^2)^{-1/2}=2.$$

In \cite{R1,R2} we proved numerous properties of local algebras.
The introduction of quasi-local algebras in \cite{GGISV} led us to
review the validity of many of these results for quasi-local
algebras. Most of the summability results of \cite{R2} are valid
in the quasi-local setting.  In addition, the summability results
of \cite{R2} are also valid for general semifinite spectral
triples since they rely only on properties of the ideals
$\LL^{(p,\infty)}$, $p\geq 1$, \cite{C,CPS2}, and the trace
property. We quote the version of the summability results from
\cite{R2} that we require below.

\begin{prop}[\cite{R2}]\label{wellbehaved} Let $(\A,\HH,\D)$ be a $QC^\infty$, local
$(1,\infty)$-summable semifinite spectral triple relative to
$(\cn,\tau)$. Let $T\in\cn$ satisfy $T\phi=\phi T=T$ for some
$\phi\in\A_c$. Then \ben T(1+\D^2)^{-1/2}\in\LL^{(1,\infty)}.\een
For $Re(s)>1$, $T(1+\D^2)^{-s/2}$ is trace class. If the limit \be
\lim_{s\to 1/2^+}(s-1/2)\tau(T(1+\D^2)^{-s})\label{mumbo}\ee
exists, then it is equal to \ben \frac{1}{2}\bigintcross
T(1+\D^2)^{-1/2}.\een In addition, for any Dixmier trace
$\tau_\omega$, the function \ben a\mapsto
\tau_\omega(a(1+\D^2)^{-1/2})\een defines a trace on
$\A_c\subset\A$.
\end{prop}

In \cite{CPRS2}, the noncommutative geometry local index theorem
of \cite{CM} was extended to semifinite spectral triples. In the
simplest terms, the local index theorem provides a formula for the
pairing of a finitely summable spectral triple $(\A,\HH,\D)$ with
the $K$-theory of $\overline{\A}$. The precise statement that we
require is
\begin{thm}[\cite{CPRS2}] Let $(\A,\HH,\D)$ be an odd $QC^\infty$
$(1,\infty)$-summable local semifinite spectral triple, relative
to $(\cn,\tau)$. Then for $u\in\A$ unitary the pairing of $[u]\in
K_1(\overline{\A})$ with $(\A,\HH,\D)$ is given by
$$ \la
[u],(\A,\HH,\D)\ra={\rm
res}_{s=0}\tau(u[\D,u^*](1+\D^2)^{-1/2-s}).$$ In particular, the
residue on the right exists.
\end{thm}
For more information on this result  see \cite{CPS2,CPRS2,CPRS3,CM}.
\vspace{-12pt}
\section{Graph $C^*$-Algebras with Semifinite Graph Traces}\label{traces}
\vspace{-12pt}
This section considers the existence of (unbounded) traces
on graph algebras. We denote by $A^+$ the positive cone in a
$C^*$-algebra $A$, and we use extended arithmetic on $[0,\infty]$ so that
 $0\times \infty=0$. From  \cite{PhR} we take the basic definition:

\begin{defn} A trace on a $C^*$-algebra $A$ is a map $\tau:A^+\to[0,\infty]$
satisfying

1) $\tau(a+b)=\tau(a)+\tau(b)$ for all $a,b\in A^+$

2) $\tau(\lambda a)=\lambda\tau(a)$ for all $a\in A^+$ and $\lambda\geq 0$

3) $\tau(a^*a)=\tau(aa^*)$ for all $a\in A$

We say: that $\tau$ is faithful if $\tau(a^*a)=0\Rightarrow a=0$; that $\tau$
 is semifinite if $\{a\in A^+:\tau(a)<\infty\}$ is norm dense in $A^+$ (or
that $\tau$ is densely defined); that $\tau$ is lower semicontinuous if
 whenever $a=\lim_{n\to\infty}a_n$ in norm in $A^+$ we have
$\tau(a)\leq\lim\inf_{n\to\infty}\tau(a_n)$.
\end{defn}
We may extend a (semifinite) trace $\tau$  by linearity to a
linear functional on  (a dense subspace of) $A$. Observe that the domain
of definition of a densely defined trace is a two-sided  ideal
$I_\tau\subset A$.

\begin{lemma}\label{finiteonfinite} Let $E$ be a row-finite directed graph and
let
$\tau:C^*(E)\to\C$ be a semifinite trace. Then the dense subalgebra
$$ A_c:={\rm span}\{S_\mu S_\nu^*:\mu,\nu\in E^*\}$$
is contained in the domain $I_\tau$ of $\tau$.
\end{lemma}

\begin{proof} Let $v\in E^0$ be a vertex, and let $p_v\in A_c$ be
 the corresponding projection. We claim that $p_v\in I_\tau$. Choose $a\in I_\tau$ positive,
so $\tau(a)<\infty$, and with $\Vert p_v-a\Vert<1$. Since $p_v$ is
a projection, we also have $\Vert p_v-p_vap_v\Vert<1$ and
$p_vap_v\in I_\tau$, so we have $\tau(p_vap_v)<\infty$.

The subalgebra $p_vC^*(E)p_v$ has unit $p_v$, and as $\Vert
p_v-p_vap_v\Vert<1$, $p_vap_v$ is invertible. Thus there is some
$b\in p_v C^*(E)p_v$ such that $bp_vap_v=p_v$. Then, again since
the trace class elements form an ideal, we have
$\tau(p_v)<\infty$.

Now since $S_\mu S_\nu^*=p_{s(\mu)}S_\mu S_\nu^*$, it is easy to see
that every element of $A_c$ has finite trace.
\end{proof}

It is convenient to denote by $A=C^*(E)$ and
$A_c=\mbox{span}\{S_\mu S_\nu^*:\mu,\nu\in E^*\}.$

\begin{lemma}\label{necessary} Let $E$ be a row-finite directed
  graph.
\par\noindent {\bf (i)} If $C^*(E)$ has a faithful semifinite trace
then no loop
can have an exit.
\par\noindent
{\bf (ii)} If $C^* (E)$ has a gauge-invariant, semifinite, lower
semicontinuous trace $\tau$ then
$\tau \circ \Phi = \tau$ and
$$
\tau(S_\mu S_\nu^*)=\delta_{\mu,\nu}\tau(p_{r(\mu)}).
$$

\noindent In particular, $\tau$ is supported on $C^* ( \{ S_\mu
S_\mu^* : \mu \in E^* \} )$.
\end{lemma}

\begin{proof}
Suppose $E$ has a loop $L = e_1 \ldots e_n$ which has an exit. Let
$v_i = s( e_i )$ for $i=1 , \cdots , n$ so that $ r ( e_n ) =
v_1$. Without loss of generality suppose that $v_1$ emits an edge
$f$ which is  not part of $L$. If $w = r(f)$ then we have
$$
\tau (p_{v_1} ) \ge \tau ( S_{e_1} S_{e_1}^* + S_f S_f^* ) = \tau
( S_{e_1}^* S_{e_1} ) + \tau ( S_f^* S_f ) = \tau ( p_{v_2} ) +
\tau ( p_w ) .
$$

\noindent Similarly we may show that $\tau ( p_{v_i} ) \ge \tau (
p_{v_{i+1}} )$ for $i = 1 , \ldots , n-1$ and so $\tau ( p_{v_1} )
\geq \tau( p_{v_1} ) + \tau ( p_w )$ which means, by Lemma
\ref{finiteonfinite}, that we must have
$\tau(p_w) =0$. Since $p_w$ is positive, this implies that $\tau$
is not faithful. Now suppose the trace $\tau$ is gauge-invariant.
Then
$$
\tau ( S_\mu S_\nu^* ) = \tau ( \gamma_z S_\mu S_\nu^* ) = \tau (
z^{\vert \mu \vert - \vert \nu \vert} S_\mu S_\nu^* ) = z^{\vert
  \mu \vert - \vert \nu \vert} \tau ( S_\mu S_\nu^* )
$$

\noindent for all $z \in S^1$, and so $\tau ( S_\mu S_\nu^* )$ is
zero unless $\vert \mu \vert = \vert \nu \vert$. Hence $\tau \circ
\Phi = \tau$ on $A_c$.
Moreover, if $\vert \mu \vert = \vert \nu \vert$ then
$$
\tau ( S_\mu S_\nu^* ) = \tau ( S_\nu^* S_\mu ) = \tau (
\delta_{\mu ,
  \nu}  p_{r ( \mu )} ) = \delta_{\mu , \nu} \tau ( p_{r ( \mu )} ) ,
$$
\noindent so  the restriction of $\tau$ to $A_c$  is supported on
 $\mbox{span} \{ S_\mu
S_\mu^* : \mu \in E^* \}$. To extend these conclusions to the
$C^*$ completions, let $\{\phi_n\}\subset \Phi(A)$ be an
approximate unit for $A$ consisting of an increasing sequence of
projections. Then for each $n$, the restriction of $\tau$ to
$A_n:=\phi_nA\phi_n$ is a finite trace, and so norm continuous.
Observe also that $\phi_nA_c\phi_n$ is dense in $A_n$ and
$\phi_nA_c\phi_n\subseteq A_c$. We claim that
\begin{equation} \mbox{when
restricted to}\ A_n,\  \tau\ \mbox{satisfies}\
\tau\circ\Phi=\tau.\label{nthoftheway}\end{equation}
 To see this we make two observations, namely that
$$\Phi(A_n)=\Phi(\phi_nA\phi_n)=\phi_n\Phi(A)\phi_n\subseteq\phi_nA\phi_n=A_n$$
and that on $\phi_nA_c\phi_n\subseteq A_c$ we have $\tau\circ\Phi=\tau$. The
norm continuity of $\tau$ on $A_n$ now completes the proof of the claim.
Now let $a\in A^+$, and let $a_n=a^{1/2}\phi_na^{1/2}$ so that $a_n\leq a_{n+1}
\leq\cdots\leq a$ and $\Vert a_n-a\Vert\to 0$. Then
$$\tau(a)\geq \lim\sup\tau(a_n)\geq\lim\inf\tau(a_n)\geq\tau(a),$$
the first inequality coming from the positivity of $\tau$, and the
last inequality from lower semicontinuity. Since $\tau$ is a trace and
$\phi_n^2=\phi_n$ we have
\begin{equation} \tau(a)=\lim_{n\to\infty}\tau(a_n)=
\lim_{n\to\infty}\tau(\phi_na\phi_n).\label{first}\end{equation}
Similarly, let $b_n=\Phi(a)^{1/2}\phi_n\Phi(a)^{1/2}$ so that
 $b_n\leq b_{n+1}\leq
\cdots\leq \Phi(a)$ and $\Vert b_n-\Phi(a)\Vert\to 0$. Then
\begin{equation}\tau(\Phi(a))=\lim_{n\to\infty}\tau(b_n)=
\lim_{n\to\infty}\tau(\phi_n\Phi(a)\phi_n)=
\lim_{n\to\infty}\tau(\Phi(\phi_na\phi_n)).\label{second}\end{equation}
However $\phi_na\phi_n\in A_n$ so by (\ref{nthoftheway}) we have
$(\tau\circ\Phi)(\phi_na\phi_n)=\tau(\phi_na\phi_n)$. Then by
Equations (\ref{first}) and (\ref{second}) we have $\tau(a)=(\tau\circ\Phi)(a)$
for all $a\in A^+$. By linearity this is true for all $a\in A$, so
$\tau=\tau\circ\Phi$ on all of $A$. Finally,
$$\phi_n\mbox{span}\{S_\mu S_\mu^*:\mu\in E^*\}\phi_n\subseteq
\mbox{span}\{S_\mu S_\mu^*:\mu\in E^*\},$$ so by the arguments
above $\tau$ is supported on $C^*(\{S_\mu S_\mu^*:\mu\in E^*\})$.
\end{proof}

Whilst the condition that no loop has an exit is necessary for the existence
of a faithful semifinite trace, it is
not sufficient.

One of the advantages of graph $C^*$-algebras is the ability to
use both graphical and analytical techniques. There is an analogue
of the above discussion of traces in terms of the graph.
\begin{defn}[cf.\ \cite{T}] If $E$ is a row-finite directed graph, then a
graph trace on $E$ is a function $g:E^0\to{\R}^+$ such that for
any $v\in E^0$ we have
\begin{equation} \label{tracecond}
g(v)=\sum_{s(e)=v}g(r(e)).
\end{equation}
\noindent If $g(v)\neq 0$ for all $v\in E^0$ we say that $g$ is
faithful.
\end{defn}

{\bf Remark} One can show by induction that if $g$ is a graph
trace on a directed graph with no sinks, and $n\geq 1$
\begin{equation} g(v)=\sum_{s(\mu)=v,\ |\mu|=n}g(r(\mu)).
\label{nosinksum}\end{equation}
For graphs with sinks, we must also count paths of length at most
$n$ which end on sinks. To deal with this more general case we
write \begin{equation} g(v)=\sum_{s(\mu)=v,\ |\mu|\preceq
n}g(r(\mu))\geq\sum_{s(\mu)=v,\ |\mu|=n}g(r(\mu)),\label{sinksum}\end{equation}
 where $|\mu|\preceq n$
means that $\mu$ is of length $n$ or is of length less than $n$
and terminates on a sink.

 As with traces on $C^*(E)$, it is easy
to see that a necessary condition for $E$ to have a faithful graph
trace is that no loop has an exit.

\begin{lemma}\label{infpaths} Suppose that $E$ is a row-finite directed
graph and
there exist vertices $v,w\in E^0$ with an infinite number of paths
from $v$ to $w$. Then there is no faithful graph trace on $E^0$.
\end{lemma}

\begin{proof} First suppose that there are an infinite number of
paths from $v$ to $w$ of the same length, $k$ say. Then for any
$N\in\N$ and any graph trace $g:E^0\to\R^+$
$$g(v)=\sum_{s(\mu)=v,\ |\mu|\preceq k}g(r(\mu))\geq \sum^N
g(w)=Ng(w).$$ So to assign a finite value to $g(v)$ we require
$g(w)=0$.

Thus we may suppose that there are infinitely many paths of
different length from $v$ to $w$, and without loss of generality
that all the paths have different length. Choose the shortest path
$\mu_1$ of length $k_1$, say. Then, with $E^m(v)=\{\mu\in
E^*:s(\mu)=v,\ |\mu|\preceq m\}$, we have \be g(v)=\sum_{\mu\in
E^{k_1}(v)}g(r(\mu))=g(w)+ \sum_{\mu\in E^{k_1}(v),\ r(\mu)\neq
w}g(r(\mu)).\label{stepone}\ee Observe that at least one of the
paths, call it $\mu_2$, in the rightmost sum can be extended until
it reaches $w$. Choose the shortest such extension from $r(\mu_2)$
to $w$, and denote the length by $k_2$. So \begin{align}
&\sum_{\mu\in E^{k_1}(v),\
\mu\neq\mu_1}g(r(\mu))=g(r(\mu_2))+\sum_{\mu\in E^{k_1}(v),\
\mu\neq\mu_1,\mu_2}g(r(\mu))\nno &=\sum_{\mu\in
E^{k_2}(r(\mu_2))}g(r(\mu))+\sum_{\mu\in E^{k_1}(v),\
\mu\neq\mu_1,\mu_2}g(r(\mu))\nno &=g(w)+\sum_{\mu\in
E^{k_2}(r(\mu_2)),\ \mu\neq\mu_2}g(r(\mu))+\sum_{\mu\in
E^{k_1}(v),\  \mu\neq\mu_1,\mu_2}g(r(\mu)).\end{align} So by
equation (\ref{stepone}) we have
$$g(v)=2g(w)+\ \mbox{sum}_1+\mbox{sum}_2.$$ The two sums on the right contain at
least one path which can be extended to $w$, and so chossing the
shortest,
$$g(v)=3g(w)+\ \mbox{sum}_1+\mbox{sum}_2+\mbox{sum}_3.$$
It is now clear how to proceed, and we deduce as before that for
all $N\in\N$, $g(v)\geq Ng(w)$.
\end{proof}

\begin{defn}\label{ends} Let $E$ be a row-finite directed graph.
 An {\em end} will mean a sink, a loop without exit or an infinite path
with no exits.
\end{defn}

{\bf Remark} We shall identify an end with the vertices which
comprise it. Once on an end (of any sort) the graph trace remains
constant.

\begin{cor} Suppose that $E$ is a row-finite directed graph and
there exists a vertex $v\in E^0$ with an infinite number of paths
from $v$ to an end. Then there is no faithful graph trace on
$E^0$.
\end{cor}

\begin{proof}  Because the value of the
graph trace is constant on an end $\Omega$, say $g_\Omega$, we
have, as in Lemma \ref{infpaths},
$$ g(v)\geq Ng_E$$
for all $N\in\N$. Hence there can be no faithful graph trace.
\end{proof}

Thus if a row-finite directed graph $E$ is to have a faithful
graph trace, it is necessary that no vertex connects infinitely
often to any other vertex or to an end, and that no loop has an
exit.

\begin{prop}\label{Eendsinends} Let $E$ be a row-finite directed
graph and suppose there exists $N\in\N$ such that for all vertices
$v$ and $w$ and for all ends $\Omega$,

1) the number of paths from $v$
to $w$, and

2) the number of paths from $v$ to $\Omega$

is less than or equal to $N$. If in addition the only infinite
paths in $E$ are eventually in ends, then $E$ has a faithful graph
trace.
\end{prop}

\begin{proof} First observe that our hypotheses on $E$ rule out
loops with exit, since we can define infinite paths using such
loops, but they are not ends.

Label the set of ends by $i=1,2,...$. Assign a positive number
$g_i$ to each end, and define $g(v)=g_i$ for all $v$ in the $i$-th
end. If there are infinitely many ends, choose the $g_i$ so that
$\sum_ig_i<\infty$.

For each end, choose a vertex $v_i$ on the end. For  $v\in E^0$
not on an end,  define
\begin{equation}g(v)=\sum_i\sum_{s(\mu)=v,\ r(\mu)=v_i}g_i.
\label{backwarddefn}\end{equation} Then the conditions on the
graph ensure this sum is finite. Using Equation (\ref{sinksum}),
one can check that $g:E^0\to\R^+$ is a faithful graph trace.
\end{proof}
There are many directed graphs with much more complicated
structure than those described in Proposition \ref{Eendsinends}
which possess faithful graph traces. The
difficulty  in defining a graph trace is going `forward', and this
is what prevents us giving a concise sufficiency condition.
Extending a graph trace `backward' from a given set of values can
always be handled as in Equation (\ref{backwarddefn}).

\begin{prop}\label{trace=graphtrace} Let $E$ be a row-finite directed graph.
Then there is a one-to-one correspondence between faithful graph
traces on $E$ and faithful, semifinite, lower semicontinuous,
 gauge invariant traces on $C^*(E)$.
\end{prop}

\begin{proof} Given a faithful graph trace $g$ on $E$ we define
$\tau_g$ on $A_c$ by
\begin{equation} \label{taudef}
\tau_g( S_\mu S_\nu^* ) %:= \tau_g (s_\nu^* s_\mu)
% :=\delta_{\mu, \nu} \tau_g ( p_{r(\mu)} )
:=\delta_{\mu , \nu} g (r( \mu ) ).
\end{equation}
One checks that $\tau_g$ is a gauge invariant trace on $A_c$,
and is faithful because for
$a = \sum_{i=1}^n c_{\mu_i , \nu_i } S_{\mu_i} S_{\nu_i}^* \in
A_c$ we have  $a^* a \ge \sum_{i=1}^n \vert c_{\mu_i , \nu_i}
\vert^2 S_{\nu_i} S_{\nu_i}^*$ and then
\begin{align}
\langle a , a \rangle_g &:= \tau_g ( a^* a ) \ge \tau_g ( \sum_{i=1}^n
\vert c_{\mu_i , \nu_i} \vert^2 S_{\nu_i} S_{\nu_i}^* ) \nno &=
\sum_{i=1}^n \vert c_{\mu_i , \nu_i} \vert^2 \tau_g ( S_{\nu_i}
S_{\nu_i}^* ) = \sum_{i=1}^n \vert c_{\mu_i , \nu_i } \vert^2 g (
r ( \nu_i ) ) > 0 .
\end{align}
Then $\la a,b\ra_g=\tau_g(b^*a)$ defines a positive definite inner product
on $A_c$ which makes it a Hilbert algebra (that the left regular
representation of $A_c$ is nondegenerate follows from $A_c^2=A_c$).

Let $\HH_g$ be the Hilbert space completion of $A_c$. Then defining
$\pi:A_c\to\B(\HH_g)$  by $\pi(a)b=ab$ for $a,b\in A_c$ yields a
faithful $*$-representation. Thus $\{\pi(S_e),\pi(p_v):e\in E^1,\ v\in E^0\}$
is a Cuntz-Krieger $E$ family in $\B(\HH_g)$. The gauge invariance
of $\tau_g$ shows that for each $z\in S^1$ the map $\gamma_z:A_c\to A_c$
 extends to a unitary $U_z:\HH_g\to\HH_g$. Then for $a,b\in A_c$ we compute
$$ (U_z\pi(a)U_{\bar{z}})(b)=U_za\gamma_{\bar{z}}(b)=
\gamma_z(a\gamma_{\bar{z}}(b))=\gamma_z(a)b=\pi(\gamma_z(a))(b).$$
Hence $U_z\pi(a)U_{\bar{z}}=\pi(\gamma_z(a))$ and defining
$\al_z(\pi(a)):=U_z\pi(a)U_{\bar{z}}$ gives a point norm
continuous action of $S^1$ on $\pi(A_c)$ implementing the gauge
action. Since for all $v\in E^0$, $\pi(p_v)p_v=p_v$, $\pi(p_v)\neq
0$. Thus we can invoke the gauge invariant uniqueness theorem,
\cite[Theorem 2.1]{BPRS}, and the map $\pi:A_c\to\B(\HH_g)$
extends by continuity to $\pi:C^*(E)\to\B(\HH_g)$ and
$\pi(C^*(E))=\overline{\pi(A_c)}^{\Vert\cdot\Vert}$ in
$\B(\HH_g)$. In particular the representation is faithful on
$C^*(E)$.

Now, $\pi(C^*(E))\subseteq\pi(A_c)''=\overline{\pi(A_c)}^{u.w.}$,
where $u.w.$ denotes the ultra-weak closure. The general theory of
Hilbert algebras, see for example  \cite[Thm 1, Sec 2, Chap 6,
Part I]{Dix}, now shows that the trace $\tau_g$ extends to an
ultra weakly lower semicontinuous, faithful, (ultra weakly)
semifinite trace $\bar{\tau}_g$ on $\pi(A_c)''$. Trivially, the
restriction of this extension to $\pi(C^*(E))$ is faithful. It is
semifinite in the norm sense on $C^*(E)$ since $\pi(A_c)$ is norm
dense in $\pi(C^*(E))$ and $\tau_g$ is finite on $\pi(A_c)$. To
see that this last statement is true, let $a\in A_c$, choose any
local unit $\phi\in A_c$ for $a$ and then
$$ \infty>\tau_g(a)=\tau_g(\phi a)=\la a,\phi\ra_g=:\bar{\tau}_g(\phi a)
=\bar{\tau}_g(a).$$
It is norm lower semicontinuous on $\pi(C^*(E))$ because if
$\pi(a)\in C^*(E)^+$ and $\pi(a_n)\in C^*(E)^+$ with $\pi(a_n)\to\pi(a)$
in norm, then $\pi(a_n)\to\pi(a)$ ultra weakly and so
$\bar{\tau}_g(\pi(a))\leq\lim\inf\bar{\tau}_g(\pi(a_n))$.

We have seen that the gauge action of $S^1$ on $C^*(E)$ is
implemented in the representation $\pi$ by the unitary
representation $S^1\ni z\to U_z\in\B(\HH_g)$. We wish to show that
$\bar{\tau}_g$ is invariant under this action, but since the $U_z$
do not lie in $\pi(A_c)''$, we can not use the tracial property
directly. Now $T\in\pi(A_c)''$ is in the domain of definition of
$\bar{\tau}_g$ if and only if $T=\pi(\xi)\pi(\eta)^*$ for left
bounded elements $\xi,\eta\in\HH_g$. Then
$\bar{\tau}_g(T)=\bar{\tau}_g(\pi(\xi)\pi(\eta)^*):=\la\xi,\eta\ra_g.$
Since $U_z(\xi)$ and $U_z(\eta)$ are also left bounded elements of
$\HH_g$ we have \bean \bar{\tau}_g(U_zTU_{\bar{z}})&=&
\bar{\tau}_g(U_z\pi(\xi)\pi(\eta)^*U_{\bar{z}})
=\bar{\tau}_g(U_z\pi(\xi)[U_z\pi(\eta)]^*)\nno
&=&\bar{\tau}_g(\pi(\gamma_z(\xi))[\pi(\gamma_z(\eta))]^*) =\la
U_z(\xi),U_z(\eta)\ra_g\nno&=&\la
\xi,\eta\ra_g=\bar{\tau}_g(T).\eean That is,
$\bar{\tau}_g(\al_z(T))=\bar{\tau}_g(T)$, and $\bar{\tau}_g$ is
$\al_z$-invariant. Thus $a\to \bar{\tau}_g(\pi(a))$ defines a
faithful, semifinite, lower semicontinuous, gauge invariant trace
on $C^*(E)$.

Conversely, given a faithful, semifinite, lower semicontinuous and gauge
invariant trace $\tau$ on $C^*(E)$, we know by Lemma \ref{finiteonfinite}
that $\tau$ is finite on $A_c$ and so we define $g(v):=\tau(p_v)$. It is easy to
 check that this is a faithful graph trace.
\end{proof}
\vspace{-24pt}
\section{Constructing a $C^*$- and Kasparov Module}\label{triplesI}
\vspace{-8pt} There are several steps in the construction of a
spectral triple. We begin in Subsection \ref{Cstarmodule} by
constructing a $C^*$-module. We define an unbounded operator $\D$
on this $C^*$-module as the generator of the gauge action of $S^1$
on the graph algebra. We show in Subsection \ref{CstarDeeee} that
$\D$ is a regular self-adjoint operator on the $C^*$-module. We
use the phase of $\D$ to construct a Kasparov module.
% we show that $\D$ is
%affiliated to a suitable semifinite von Neumann algebra with a
%faithful normal trace, and has compact resolvent in the von
%Neumann sense.

\subsection{Building a $C^*$-module}\label{Cstarmodule}
\vspace{-7pt} The constructions of this subsection work for any
locally finite graph.  Let $A=C^*(E)$ where $E$ is any locally
finite directed graph. Let $F=C^*(E)^\gamma$ be the fixed point
subalgebra for the gauge action. Finally, let $A_c,F_c$ be the
dense subalgebras of $A,F$ given by the (finite) linear span of
the generators.

We make $A$ a right inner product $F$-module. The right action of
$F$ on $A$ is by right multiplication. The inner product is
defined by
$$ (x|y)_R:=\Phi(x^*y)\in F.$$
Here $\Phi$ is the canonical expectation. It is simple to check
the requirements that $(\cdot|\cdot)_R$ defines an $F$-valued
inner product on $A$. The requirement $(x|x)_R=0\Rightarrow x=0$
follows from the faithfulness of $\Phi$.

\begin{defn}\label{Fmod} Define $X$ to be the $C^*$-$F$-module
completion of $A$
for the $C^*$-module norm
$$\Vert x\Vert_X^2:=\Vert(x|x)_R\Vert_A=\Vert(x|x)_R\Vert_F=
\Vert \Phi(x^*x)\Vert_F.$$
Define $X_c$ to be the pre-$C^*$-$F_c$-module with linear space
$A_c$ and the inner product $(\cdot|\cdot)_R$.
\end{defn}

{\bf Remark} Typically, the action of $F$ does not map $X_c$ to
itself, so we may only consider $X_c$ as an $F_c$ module. This is
a reflection of the fact that $F_c$ and $A_c$ are quasilocal not
local.

The inclusion map $\iota:A\to X$ is continuous since
$$\Vert a\Vert_X^2=\Vert\Phi(a^*a)\Vert_F\leq\Vert
a^*a\Vert_A=\Vert a\Vert^2_A.$$ We can also define the gauge
action $\gamma$ on $A\subset X$, and as
\bean\Vert\gamma_z(a)\Vert^2_X&=&\Vert\Phi((\gamma_z(a))^*(\gamma_z(a)))\Vert_F
=\Vert\Phi(\gamma_z(a^*)\gamma_z(a))\Vert_F\nno&=&
\Vert\Phi(\gamma_z(a^*a))\Vert_F =\Vert\Phi(a^*a)\Vert_F=\Vert
a\Vert^2_X,\eean for each $z\in S^1$, the action of $\gamma_z$ is
isometric on $A\subset X$ and so extends to a unitary $U_z$ on
$X$. This unitary is $F$ linear, adjointable, and we obtain a
strongly continuous action of $S^1$ on $X$, which we still denote
by $\gamma$.

For each $k\in\Z$, the projection onto the $k$-th spectral subspace for
the gauge action defines an operator $\Phi_k$ on $X$  by
$$\Phi_k(x)=
\frac{1}{2\pi}\int_{S^1}z^{-k}\gamma_z(x)d\theta,\ \ z=e^{i\theta},\ \ x\in X.$$
Observe that on generators we have $\Phi_k(S_\al
S_\beta^*)=S_\al S_\beta^*$ when
$|\al|-|\beta|=k$ and is zero when $|\al|-|\beta|\neq k$. The range of $\Phi_k$
is
\begin{equation}
\mbox{Range}\ \Phi_k=\{x\in X:\gamma_z(x)=z^kx\ \ \mbox{for all}\
z\in S^1\}. \label{kthproj}\end{equation} These ranges give us a
natural $\Z$-grading of $X$.

{\bf Remark} If $E$ is a finite graph with no loops, then for $k$
sufficiently large there are no paths of length $k$ and so
$\Phi_k=0$. This will obviously simplify many of the convergence
issues below.

\begin{lemma}\label{phiendo} The operators $\Phi_k$ are adjointable endomorphisms
of the $F$-module $X$ such that $\Phi_k^*=\Phi_k=\Phi_k^2$ and
$\Phi_k\Phi_l=\delta_{k,l}\Phi_k$. If $K\subset\Z$ then the sum
$\sum_{k\in K}\Phi_k$ converges strictly to a projection in the
endomorphism algebra. The sum $\sum_{k\in\Z}\Phi_k$ converges to
the identity operator on $X$.
\end{lemma}

\begin{proof} It is clear from the definition that each $\Phi_k$ defines an
$F$-linear map on $X$. First, we show that $\Phi_k$ is bounded:
$$\Vert\Phi_k(x)\Vert_X\leq\frac{1}{2\pi}
\int_{S^1}\Vert\gamma_z(x)\Vert_Xd\theta\leq\frac{1}{2\pi}\int_{S^1}\Vert
x\Vert_Xd\theta=\Vert x\Vert_X.$$ So $\Vert \Phi_k\Vert\leq 1$.
Since $\Phi_k S_\mu=S_\mu$ whenever $\mu$ is a path of length $k$,
$\Vert\Phi_k\Vert=1$.

%To see that for $k\neq l$ $\Phi_k$ and $\Phi_l$ commute, we use
%the continuity of $S^1\times S^1\ni(z,w)\to
%z^{-l}w^{-k}\gamma_z(x)\gamma_w(x)$ for fixed $x\in X$. Since
%continuity implies measurability, we can apply Fubini to deduce
%that \bean \Phi_k\Phi_l(x)&=&\frac{1}{(2\pi)^2}
%\int_{S^1}w^{-k}\gamma_w\left(z^{-l}\gamma_z(x)dz\right)dw
%=\frac{1}{(2\pi)^2}\int_{S^1}\int_{S^1}w^{-k}z^{-l}\gamma_{wz}(x)dzdw\nno
%&=&\frac{1}{(2\pi)^2}\int_{S^1}\int_{S^1}z^{-l}w^{-k}\gamma_{zw}(x)dwdz
%=\Phi_l\Phi_k(x).\eean
On the subspace $X_c$ of finite linear combinations of generators,
one can use Equation (\ref{kthproj}) to see that
$\Phi_k\Phi_l=\delta_{k,l}\Phi_k$ since
$$\Phi_k\Phi_lS_\al S_\beta^*=\Phi_k\delta_{|\al|-|\beta|,l}S_\al
S_\beta^*=\delta_{|\al|-|\beta|,k}\delta_{|\al|-|\beta|,l}S_\al
S_\beta^*.$$ For general $x\in X$, we approximate $x$ by a
sequence $\{x_m\}\subset X_c$, and the continuity of the $\Phi_k$
then shows that the relation $\Phi_k\Phi_l=\delta_{k,l}\Phi_k$
holds on all of $X$. Again using the continuity of $\Phi_k$, the
following computation allows us to show that for all $k$, $\Phi_k$
is adjointable  with adjoint $\Phi_k$: \bean (\Phi_kS_\al
S_\beta^*|S_\rho S_\s^*)_R&=&
\Phi\left(\delta_{|\al|-|\beta|,k}S_\beta S_\al^*S_\rho
S_\s^*\right)\nno
&=&\delta_{|\al|-|\beta|,k}\delta_{|\beta|-|\al|+|\rho|-|\s|,0}S_\beta
S_\al^*S_\rho S_\s^*\nno
&=&\Phi\left(\delta_{|\rho|-|\s|,k}S_\beta S_\al^*S_\rho
S_\s^*\right) =(S_\al S_\beta^*|\Phi_{k}S_\rho S_\s^*)_R.\eean

To address the last two statements of the Lemma, we observe that
the set $\{\Phi_k\}_{k\in\Z}$ is norm bounded in $End_F(X)$, so
the strict topology on this set coincides with the $*$-strong
topology, \cite[Lemma C.6]{RW}. First, if $K\subset\Z$ is a finite
set, the sum
$$\sum_{k\in K}\Phi_k$$
is finite, and defines a projection in $End_F(X)$ by the results
above. So assume $K$ is infinite and let $\{K_i\}$ be an
increasing sequence of finite subsets of $K$ with $K=\cup_iK_i$.
For $x\in X$, let
$$T_ix=\sum_{k\in K_i}\Phi_kx.$$
Choose a sequence $\{x_m\}\subset X_c$ with $x_m\to x$. Let
$\epsilon>0$ and choose $m$ so that $\Vert
x_m-x\Vert_X<\epsilon/2$. Since $x_m$ has finite support, for
$i,j$ sufficiently large we have $T_ix_m-T_jx_m=0$, and so for
sufficiently large $i,j$ \bean \Vert T_ix-T_jx\Vert_X&=&\Vert
T_ix-T_ix_m+T_ix_m-T_jx_m+T_jx_m-T_jx_m\Vert_X\nno &\leq&\Vert
T_i(x-x_m)\Vert_X+\Vert T_j(x-x_m)\Vert_X+\Vert
T_ix_m-T_jx_m\Vert_X\nno &<&\epsilon.\eean This proves the strict
convergence, since the $\Phi_k$ are all self-adjoint. To prove the
final statement, let $x,\{x_m\}$ be as above, $\epsilon>0$, and
choose $m$ so that $\Vert x-x_m\Vert_X<\epsilon/2$. Then \bean
\Vert x-\sum_{k\in\Z}\Phi_kx\Vert_X&=&\Vert x-\sum\Phi_k
x_m+\sum\Phi_kx_m-\sum\Phi_kx\Vert_X\nno &\leq&\Vert
x-x_m\Vert_X+\Vert\sum\Phi_k(x-x_m)\Vert_X<\epsilon.\qed\eean
\hideqed
\end{proof}

\begin{cor}\label{gradedsum} Let $x\in X$. Then with $x_k=\Phi_kx$ the sum
$\sum_{k\in \Z}x_k$ converges in $X$ to $x$.
\end{cor}
\vspace{-7pt}
\subsection{The Kasparov Module}\label{CstarDeeee}
\vspace{-7pt} {\bf In this subsection we assume that $E$ is
locally finite and furthermore has no sources. That is, every
vertex receives at least one edge.}

Since we have the gauge action defined on $X$, we may use the
generator of this action to define an unbounded operator $\D$. We
will not define or study $\D$ from the generator point of view,
rather taking a more bare-hands approach. It is easy to check that
$\D$ as defined below is the generator of the $S^1$ action.

 The theory of unbounded operators on
$C^*$-modules that we require is all contained in Lance's book,
\cite[Chapters 9,10]{L}. We quote the following definitions
(adapted to our situation).
\begin{defn} Let $Y$ be a right $C^*$-$B$-module. A densely defined
 unbounded operator $\D:{\rm dom}\ \D\subset Y\to Y$ is a
 $B$-linear operator defined on a dense $B$-submodule
 ${\rm dom}\ \D\subset Y$. The operator $\D$ is closed
 if the graph
 $$ G(\D)=\{(x|\D x)_R:x\in{\rm dom}\ \D\}$$
 is a closed submodule of $Y\oplus Y$.
\end{defn}
If $\D:\mbox{dom}\ \D\subset Y\to Y$ is densely defined and
 unbounded, define
 a submodule
 $$\mbox{dom}\ \D^*:=\{y\in Y:\exists z\in Y\ \mbox{such that}\
 \forall x\in\mbox{dom}\ \D, (\D x|y)_R=(x|z)_R\}.$$
 Then for $y\in \mbox{dom}\ \D^*$ define $\D^*y=z$.
Given $y\in\mbox{dom}\ \D^*$,
 the element $z$ is unique, so $\D^*:\mbox{dom}\D^*\to Y$, $\D^*y=z$
is well-defined, and
 moreover is closed.
 \begin{defn} Let $Y$ be a right $C^*$-$B$-module. A densely defined unbounded
 operator $\D:{\rm dom}\ \D\subset Y\to Y$ is symmetric if for all
 $x,y\in{\rm dom}\ \D$
 $$ (\D x|y)_R=(x|\D y)_R.$$
 A symmetric operator $\D$ is self-adjoint if
 ${\rm dom}\ \D={\rm dom}\ \D^*$ (and so $\D$ is necessarily
 closed). A densely defined unbounded operator $\D$ is regular if
 $\D$ is closed, $\D^*$ is densely defined, and $(1+\D^*\D)$ has
 dense range.
 \end{defn}
 The extra requirement of regularity is necessary in the
 $C^*$-module context for the continuous functional calculus,
 and is not automatic,
 \cite[Chapter 9]{L}.

With these definitions in hand, we return to our $C^*$-module $X$.
\begin{prop}\label{CstarDee} Let $X$ be the right $C^*$-$F$-module of
Definition
 \ref{Fmod}.  Define $X_\D\subset X$ to be the linear space
$$ X_\D=
\{x=\sum_{k\in\Z}x_k\in X:\Vert\sum_{k\in\Z}k^2(x_k|x_k)_R\Vert<\infty\}.$$
For $x=\sum_{k\in\Z}x_k\in X_\D$ define
$$ \D x=\sum_{k\in\Z}kx_k.$$
Then $\D:X_\D\to X$ is a self-adjoint regular operator on $X$.
\end{prop}

{\bf Remark} Any  $S_\al S_\beta^*\in A_c$ is in $X_\D$ and
$$\D S_\al S_\beta^*=(|\al|-|\beta|)S_\al S_\beta^*.$$

\begin{proof} First we show that $X_\D$ is a submodule. If $x\in X_\D$
and $f\in F$, in the $C^*$-algebra $F$  we have
\bean\sum_{k\in\Z}k^2(x_kf|x_kf)_R&=&\sum_{k\in\Z}k^2f^*(x_k|x_k)_Rf
=f^*\sum_{k\in\Z}k^2(x_k|x_k)_Rf\nno &\leq&
 f^*f\Vert\sum_{k\in\Z}k^2 (x_k|x_k)_R\Vert.\eean
So
$$\Vert \sum_{k\in\Z}k^2(x_kf|x_kf)_R\Vert
\leq\Vert
f^*f\Vert\ \Vert\sum_{k\in\Z}k^2(x_k|x_k)_R\Vert<\infty.$$
Observe that if $x\in X$ is a finite sum of graded components,
$$ x=\sum_{k=-N}^Mx_k,$$
then $x\in X_\D$. In particular if $P=\sum_{finite}\Phi_k$ is a finite
sum of the projections $\Phi_k$, $Px\in X_\D$ for any $x\in X$.

The following calculation shows that $\D$ is symmetric on its
domain, so that the adjoint is densely defined. Let
$x,y\in\mbox{dom}\D$ and use Corollary \ref{gradedsum} to write
$x=\sum_kx_k$ and $y=\sum_ky_k$. Then \bean (\D
x|y)_R&=&(\sum_kkx_k|\sum_my_m)_R=\Phi((\sum_kkx_k)^*(\sum_my_m))
=\Phi(\sum_{k,m}kx_k^*y_m)\nno&=&\sum_kkx_k^*y_k
=\Phi(\sum_{k,m}x_m^*ky_k)=\Phi((\sum_mx_m)^*(\sum_kky_k))\nno
&=&(x|\D y)_R.\eean Thus $\mbox{dom}\D\subseteq\mbox{dom}\D^*$, and so
$\D^*$ is densely defined, and  of course closed.
Now choose any  $x\in X$ and any $y\in \mbox{dom}\D^*$. Let
$P_{N,M}=\sum_{k=-N}^M\Phi_k$, and recall that $P_{N,M}x\in\mbox{dom}\D$
for all $x\in X$. Then
\bean (x|P_{N,M}\D^*y)_R=(P_{N,M}x|\D^*y)_R&=&(\D P_{N,M}x|y)_R\nno
&=&(\sum_{k=-N}^Mkx_k|y)_R=(x|\sum_{k=-N}^Mky_k)_R.\eean
Since this is true for all $x\in X$ we have
$$P_{N,M}\D^*y=\sum_{k=-N}^Mky_k.$$
Letting $N,M\to\infty$, the limit on the left hand side exists by
Corollary \ref{gradedsum},
and so the
limit on the right exists, and so $y\in\mbox{dom}\D$. Hence $\D$ is
self-adjoint.

Finally, we need to show that $\D$ is regular. By \cite[Lemma
9.8]{L}, $\D$ is regular if and only if the operators $\D\pm
iId_X$ are surjective. This is straightforward though, for if
$x=\sum_kx_k$ we have
$$ x=\sum_{k\in\Z} \frac{(k\pm i)}{(k\pm i)}x_k=
(\D\pm iId_X)\sum_{k\in\Z}\frac{1}{(k\pm
i)}x_k.$$ The convergence of $\sum_kx_k$ ensures the convergence
of $\sum_k(k\pm i)^{-1}x_k$.
\end{proof}

There is a continuous functional calculus for self-adjoint regular
operators, \cite[Theorem 10.9]{L}, and we use this to obtain
spectral projections for $\D$ at the $C^*$-module level. Let
$f_k\in C_c({\R})$ be $1$ in a small neighbourhood of $k\in{\Z}$
and zero on $(-\infty,k-1/2]\cup[k+1/2,\infty)$. Then it is clear
that
$$ \Phi_k=f_k(\D).$$
That is the spectral projections of $\D$ are the same as the
projections onto the spectral subspaces of the gauge action.

The next Lemma is the first place where we need our graph to be
locally finite and have no sources.
\begin{lemma}\label{finrank} Assume that the directed graph $E$ is
locally finite and has no sources.
For all $a\in A$ and $k\in\Z$, $a\Phi_k\in End^0_F(X)$, the
 compact  endomorphisms of the right $F$-module $X$. If $a\in A_c$ then
$a\Phi_k$ is finite rank.
\end{lemma}

{\bf Remark} The proof actually shows that for $k>0$
$$\Phi_k=\sum_{|\rho|=k}\Theta^R_{S_\rho,S_\rho}$$
where the sum converges in the strict topology.

\begin{proof} We will prove the Lemma by first showing that for each $v\in E^0$
and $k\geq 0$
$$p_v\Phi_k=\sum_{s(\rho)=v,\ |\rho|=k}\Theta^R_{S_\rho,S_\rho}.$$
This is a finite sum, by the row-finiteness of $E$. For $k<0$ the situation is
more complicated, but a similar formula holds in that case also.

First suppose that $k\geq 0$ and $a=p_v\in A_c$ is the projection
corresponding to a vertex $v\in E^0$.  For $\al$ with $|\al|\geq
k$ denote by $\underline{\al}=\al_1\cdots\al_{k}$ and
$\overline{\al}=\al_{k+1}\cdots\al_{|\al|}$. With this notation we
compute the action of $p_v$ times the  rank one endomorphism
$\Theta^R_{S_\rho,S_\rho}$, $|\rho|=k$, on $S_\al S_\beta^*$. We
find \bean p_v\Theta^R_{S_\rho,S_\rho}S_\al
S_\beta^*&=&p_vS_\rho(S_\rho|S_\al
S_\beta^*)_R=\delta_{v,s(\rho)}p_vS_\rho\Phi(S_\rho^*S_\al
S_\beta^*)\nno
&=&\delta_{v,s(\rho)}p_vS_\rho\delta_{|\al|-|\beta|,k}
\delta_{\rho,\underline{\al}} S_{\overline{\al}}S_\beta^*=
\delta_{|\al|-|\beta|,k}\delta_{\rho,\underline{\al}}\delta_{v,s(\rho)}S_\al
S_\beta^*.\eean Of course if $|\al|<|\rho|$ we have
$$p_v\Theta^R_{S_\rho,S_\rho}S_\al S_\beta^*=
p_vS_\rho\Phi(S_\rho^*S_\al S_\beta^*)=0.$$
This too is $\delta_{|\al|-|\beta|,k}p_vS_\al S_\beta^*$.
Thus for any $\al$ we have
$$ \sum_{|\rho|=k}p_v\Theta^R_{S_\rho,S_\rho}S_\al
S_\beta^*= \sum_{|\rho|=k,
s(\rho)=v}\delta_{v,s(\rho)}\delta_{|\al|-|\beta|,k}
\delta_{\rho,\underline{\al}}p_vS_\al
S_\beta^*=\delta_{v,s(\al)}\delta_{|\al|-|\beta|,k}S_\al
S_\beta^*.$$ This is of course the action of $p_v\Phi_k$ on $S_\al
S_\beta^*$, and if $v$ is a sink, $p_v\Phi_k=0$, as it must. Since
$E$ is locally finite, the number of paths of length $k$ starting
at $v$ is finite, and we have a finite sum. For general $a\in A_c$
we may write
$$a=\sum_{i=1}^nc_{\mu_i,\nu_i}S_{\mu_i}S^*_{\nu_i}$$ for some
paths $\mu_i,\nu_i$. Then
$S_{\mu_i}S^*_{\nu_i}=S_{\mu_i}S^*_{\nu_i}p_{s(\nu_i)}$, and we
may apply the above reasoning to each term in the sum defining $a$
to get a finite sum again. Thus $a\Phi_k$ is finite rank.

Now we consider $k<0$. Given $v\in E^0$, let $|v|_k$ denote the
number of paths $\rho$ of length $|k|$ ending at $v$, i.e.
$r(\rho)=v$. Since  we assume that $E$ is locally finite and has
no sources, $\infty>|v|_k>0$ for each $v\in E^0$. We consider the
action of the finite rank operator
$$
\frac{1}{|v|_k}\sum_{|\rho|=|k|,r(\rho)=v}p_v\Theta^R_{S^*_\rho,S^*_\rho}.$$
For $S_\al S_\beta^*\in X$ we find \bean
\frac{1}{|v|_k}\sum_{|\rho|=|k|,r(\rho)=v}p_v\Theta^R_{S^*_\rho,S^*_\rho}S_\al
S_\beta^*&=&\frac{1}{|v|_k}\sum_{|\rho|=|k|,r(\rho)=v}p_vS_\rho^*\Phi(S_\rho
S_\al S_\beta^*)\nno &=&\frac{1}{|v|_k}\sum_{|\rho|=|k|,r(\rho)=v}
\delta_{|\al|-|\beta|,-|k|}p_vS_\rho^*S_\rho S_\al S_\beta^*\nno
&=&\delta_{|\al|-|\beta|,-|k|}\delta_{v,s(\al)}p_v S_\al
S_\beta^*=p_v\Phi_{k}S_\al S_\beta^*.\eean Thus $p_v\Phi_{-|k|}$
is a finite rank endomorphism, and by the argument above, we have
$a\Phi_{-|k|}$ finite rank for all $a\in A_c$.
To see that $a\Phi_k$ is compact for all $a\in A$, recall that
every $a\in A$ is a norm limit of a sequence $\{a_i\}_{i\geq
0}\subset A_c$. Thus for any $k\in\Z$
$a\Phi_k=\lim_{i\to\infty}a_i\Phi_k$ and so is compact.
\end{proof}

\begin{lemma}\label{compactendo} Let $E$ be a locally finite directed graph
with no sources. For all $a\in A$,
$a(1+\D^2)^{-1/2}$ is a compact
 endomorphism of the $F$-module $X$.
\end{lemma}

\begin{proof}
First let $a=p_v$ for $v\in E^0$. Then the sum
$$ R_{v,N}:=p_v\sum_{k=-N}^N\Phi_k(1+k^2)^{-1/2}$$
is finite rank, by Lemma \ref{finrank}. We will show that the
sequence $\{R_{v,N}\}_{N\geq 0}$  is convergent with respect to
the operator norm $\Vert\cdot\Vert_{End}$ of endomorphisms of $X$.
Indeed, assuming that $M>N$, \begin{align} \Vert
R_{v,N}-R_{v,M}\Vert_{End}&=\Vert
p_v\sum_{k=-M}^{-N}\Phi_k(1+k^2)^{-1/2}+p_v\sum_{k=N}^M\Phi_k(1+k^2)^{-1/2}\Vert_{End}\nno
&\leq2(1+N^2)^{-1/2}\to 0,\end{align} since the ranges of the
$p_v\Phi_k$ are orthogonal for different $k$. Thus, using the
argument from Lemma \ref{finrank}, $a(1+\D^2)^{-1/2}\in
End^0_F(X)$. Letting $\{a_i\}$ be a Cauchy sequence  from $A_c$,
we have
$$\Vert a_i(1+\D^2)^{-1/2}-a_j(1+\D^2)^{-1/2}\Vert_{End}\leq\Vert
a_i-a_j\Vert_{End}=\Vert a_i-a_j\Vert_A\to 0,$$ since
$\Vert(1+\D^2)^{-1/2}\Vert\leq 1$. Thus the sequence
$a_i(1+\D^2)^{-1/2}$ is Cauchy in norm and
 we
see that $a(1+\D^2)^{-1/2}$ is compact for all $a\in A$.
\end{proof}

\begin{prop}\label{Kasmodule} Assume that the directed graph $E$ is
locally finite and has no sources. Let $V=\D(1+\D^2)^{-1/2}$. Then
$(X,V)$ defines a class in $KK^1(A,F)$.
\end{prop}

\begin{proof} We will use the approach of \cite[Section 4]{K}.
We need to show that various operators belong to $End^0_F(X)$.
First, $V-V^*=0$, so $a(V-V^*)$ is compact for all $a\in A$. Also
$a(1-V^2)=a(1+\D^2)^{-1}$ which is compact from Lemma
\ref{compactendo} and the boundedness of $(1+\D^2)^{-1/2}$.
Finally, we need to show that $[V,a]$ is compact for all $a\in A$.
First we suppose that $a\in A_c$. Then
 \bean
[V,a]&=&[\D,a](1+\D^2)^{-1/2}-\D(1+\D^2)^{-1/2}[(1+\D^2)^{1/2},a](1+\D^2)^{-1/2}\nno
&=&b_1(1+\D^2)^{-1/2}+Vb_2(1+\D^2)^{-1/2},\eean where
$b_1=[\D,a]\in A_c$ and $b_2=[(1+\D^2)^{1/2},a]$. Provided that
$b_2(1+\D^2)^{-1/2}$ is a compact endomorphism,
Lemma \ref{compactendo} will show
that $[V,a]$ is compact for all $a\in A_c$. So consider the action
of $[(1+\D^2)^{1/2},S_\mu S_\nu^*](1+\D^2)^{-1/2}$ on $x=\sum_{k\in\Z}x_k$. We
find \bea &&\sum_{k\in\Z}[(1+\D^2)^{1/2},S_\mu
S_\nu^*](1+\D^2)^{-1/2}x_k\nno
&=&
\sum_{k\in\Z}
\left((1+(|\mu|-|\nu|+k)^2)^{1/2}-(1+k^2)^{1/2}\right)(1+k^2)^{-1/2}S_\mu
S_\nu^*x_k\nno &=&\sum_{k\in\Z}f_{\mu,\nu}(k)S_\mu S_\nu^*
\Phi_kx.\label{limit}\eea The function \ben
f_{\mu,\nu}(k)=\left((1+(|\mu|-|\nu|+k)^2)^{1/2}-(1+k^2)^{1/2}\right)
(1+k^2)^{-1/2}\een
goes to $0$  as $k\to\pm\infty$, and as the $S_\mu S_\nu^*\Phi_k$
are finite rank with orthogonal ranges, the sum in (\ref{limit})
converges in the endomorphism norm, and so converges to a compact
endomorphism. For  $a\in A_c$ we write $a$ as a finite
linear combination of generators $S_\mu S_\nu^*$, and apply the
above reasoning to each term in the sum to find that
$[(1+\D^2)^{1/2},a](1+\D^2)^{-1/2}$ is a compact endomorphism.
Now let $a\in A$ be the norm limit of a Cauchy sequence
$\{a_i\}_{i\geq 0}\subset A_c$. Then
$$\Vert[V,a_i-a_j]\Vert_{End}\leq 2\Vert a_i-a_j\Vert_{End}\to 0,$$
so the sequence $[V,a_i]$ is also Cauchy in norm, and so the limit
is compact.
\end{proof}

\section{The Gauge Spectral Triple of a Graph Algebra}\label{triplesII}

In this section we will construct a semifinite spectral triple for
those graph $C^*$-algebras which possess a faithful gauge
invariant trace, $\tau$. Recall from Proposition
\ref{trace=graphtrace} that  such traces arise from
 faithful graph traces.

We will begin with the right $F_c$ module $X_c$. In order to deal
with the spectral projections of $\D$ we will also assume
throughout this section that $E$ is locally finite and has no
sources. This ensures, by Lemma \ref{finrank} that for all $a\in
A$ the endomorphisms $a\Phi_k$ of $X$ are compact endomorphisms.

As in the proof of Proposition \ref{trace=graphtrace}, we
define a ${\C}$-valued inner product on $X_c$:
$$ \la x,y\ra:=\tau((x|y)_R)=\tau(\Phi(x^*y))=\tau(x^*y).$$
This inner product is linear in the second variable. We define the
Hilbert space $\HH=L^2(X,\tau)$ to be the completion of $X_c$ for
$\la\cdot,\cdot\ra$. We need a few lemmas in order to obtain the
ingredients of our spectral triple.

%\begin{proof} This is just a calculation:
%$$\tau((x,y)_L)=\tau(\Phi(xy^*))=\tau(\Phi(yx^*)^*)
%=\overline{\tau(\Phi(yx^*))}=\overline{\tau(yx^*)}
%=\overline{\tau(x^*y)}=\overline{\tau(\Phi(x^*y))}=\overline{\la
%x,y\ra}.$$
%\end{proof}

\begin{lemma}\label{endoproof} The $C^*$-algebra $A=C^*(E)$ acts on $\HH$
by an extension of left multiplication. This defines a faithful nondegenerate
$*$-representation of $A$. Moreover, any endomorphism of $X$ leaving $X_c$
invariant extends uniquely to a bounded linear operator on $\HH$.
\end{lemma}

\begin{proof} The first statement follows from the proof of
Proposition \ref{trace=graphtrace}. Now let $T$
be an endomorphism of  $X$ leaving $X_c$ invariant. Then
\cite[Cor 2.22]{RW},
$$(Tx|Ty)_R\leq \| T\|_{End}^2(x|y)_R$$
in the algebra $F$. Now the  norm of $T$ as an operator on $\HH$,
denoted $\Vert T\Vert_\infty$, can be computed in terms of the
endomorphism norm of $T$ by \begin{align}
\|T\|_\infty^2&:=\sup_{\|x\|_\HH\leq 1}\la
Tx,Tx\ra=\sup_{\|x\|_\HH\leq 1}\tau((Tx|Tx)_R)\nno &\leq
\sup_{\|x\|_\HH\leq 1}\n T\n_{End}^2\tau((x|x)_R)=\n
T\n_{End}^2.\qed\end{align} \hideqed
\end{proof}
%\begin{cor}\label{bddendos} All  endomorphisms of $X$ extend to  bounded
%operators on $\HH$. Thus, the multiplication
%representation extends to the multiplier algebra of $A$.
%\end{cor}
%\begin{proof} The proof of Lemma \ref{endoproof} applies to any
%endomorphism of $X$.
%\end{proof}
\begin{cor} The endomorphisms $\{\Phi_k\}_{k\in\Z}$ define mutually orthogonal projections on
$\HH$. For any $K\subset \Z$ the sum $\sum_{k\in K}\Phi_k$
converges strongly to a projection in $\B(\HH)$. In particular,
$\sum_{k\in\Z}\Phi_k=Id_{\HH}$, and for all $x\in \HH$ the sum
$\sum_k\Phi_kx$ converges in norm to $x$.
\end{cor}
\begin{proof} As in Lemma \ref{phiendo}, we can use the
continuity of the $\Phi_k$ on $\HH$, which follows from Corollary
\ref{endoproof}, to see that the relation
$\Phi_k\Phi_l=\delta_{k,l}\Phi_k$ extends from $X_c\subset\HH$ to
$\HH$. The  strong convergence of sums of $\Phi_k$'s
is just as in Lemma \ref{phiendo} after  replacing the
$C^*$-module norm with the Hilbert space norm.
\end{proof}
\begin{lemma} The operator $\D$ restricted to $X_c$ extends to a
closed self-adjoint
operator on $\HH$.
\end{lemma}
\begin{proof}
%We define the domain of $\D$ to be the
%completion of $X_c$ in the norm
%$$x\to\Vert x\Vert_{\HH,\D}:=\Vert
%x\Vert_{\HH}+\Vert\D x\Vert_{\HH}.$$
The proof is essentially the same as
Proposition \ref{CstarDee}.
\end{proof}

%The Hilbert space $\HH$ and operator $\D$ are two of the
%ingredients of our spectral triple. We also need a $*$-algebra. In
%fact $A_c$ will do the job, but it also has a natural completion
%$\A$  which is useful too.

\begin{lemma}\label{deltacomms} Let $\HH,\D$ be as above and let
$\dd=\sqrt{\D^*\D}=\sqrt{\D^2}$ be the absolute value of $\D$.
Then for $S_\al S_\beta^*\in A_c$, the operator $[\dd,S_\al
S_\beta^*]$ is well-defined on $X_c$, and extends to a bounded
operator on $\HH$ with
$$\Vert[\dd,S_\al S_\beta^*]\Vert_{\infty}\leq
\Bigl||\al|-|\beta|\Bigr|.$$ Similarly, $\Vert[\D,S_\al
S_\beta^*]\Vert_\infty= \Bigl||\al|-|\beta|\Bigr|$.
\end{lemma}

\begin{proof} It is clear that $S_\al S_\beta^*X_c\subset X_c$, so
we may define the action of the commutator on elements of $X_c$.
Now let $x=\sum_kx_k\in\HH$ and consider the action of $[\dd,S_\al
S_\beta^*]$ on $x_k$. We have
$$[\dd,S_\al
S_\beta^*]x_k=\Bigl(\Bigl||\al|-|\beta|+k\Bigr|-\Bigl|k\Bigr|\Bigr)S_\al
S_\beta^*x_k,$$ and so, by the triangle inequality,
$$\Vert[\dd,S_\al
S_\beta^*]x_k\Vert_{\infty}\leq\Bigl||\al|-|\beta|\Bigr|\Vert
x_k\Vert_\infty,$$ since $\Vert S_\al S_\beta^*\Vert_\infty=1.$ As the
$x_k$ are mutually orthogonal,
 $\Vert[\dd,S_\al S_\beta^*]\Vert_\infty\leq
\Bigl||\al|-|\beta|\Bigr|$. The statements about  $[\D,S_\al
S_\beta^*]=(|\al|-|\beta|)S_\al S_\beta^*$ are easier.
\end{proof}
\begin{cor}\label{smodense} The algebra $A_c$ is contained in the smooth domain
of the derivation $\delta$ where for $T\in\B(\HH)$,
$\delta(T)=[\dd,T]$. That is
$$ A_c\subseteq\bigcap_{n\geq 0}{\rm dom}\ \delta^n.$$
\end{cor}
\begin{defn} Define the $*$-algebra $\A\subset A$ to be the
completion of $A_c$ in the $\delta$-topology. By Lemma \ref{smo},
$\A$ is Fr\'{e}chet and stable under the holomorphic functional
calculus.
\end{defn}
\begin{lemma}\label{smoalg} If $a\in\A$ then $[\D,a]\in\A$ and the operators $\delta^k(a)$,
$\delta^k([\D,a])$ are bounded for all $k\geq 0$. If $\phi\in F\subset\A$ and
$a\in\A$
satisfy $\phi a=a=a\phi$, then $\phi[\D,a]=[\D,a]=[\D,a]\phi$.
The norm closed algebra generated by $\A$ and $[\D,\A]$ is $A$. In
particular, $\A$ is quasi-local.
\end{lemma}
We leave the straightforward proofs of these statements to the
reader.

%At this point we have most of the structure of a
%semifinite local spectral triple. We require only the compactness of
%$a(\lambda-\D)^{-1}$,
%$\lambda\in\C\setminus\R$, $a\in\A$, relative to some trace $\tilde\tau$ on
%some von Neumann algebra $\cn$ to which $\D$ is affiliated.
%There is a
%canonical choice of von Neumann algebra and trace, and for this
%choice $a(1+\D^2)^{-1/2}$ is in the domain of the Dixmier trace
%for all $a\in\A$.

\vspace{-7pt}
\subsection{Traces and Compactness Criteria}
\vspace{-9pt}
 We still assume that $E$ is a locally finite
graph with no sources and that $\tau$ is a faithful semifinite lower
semicontinuous gauge
invariant trace on $C^*(E)$. We will define a von Neumann algebra $\cn$
with a faithful semifinite normal trace $\tilde\tau$ so that
$\A\subset\cn\subset\B(\HH)$, where $\A$ and $\HH$ are as defined in the last
subsection. Moreover the operator $\D$ will be affiliated to $\cn$.
The aim of this subsection will then be  to prove the following result.
\begin{thm}\label{mainthm} Let $E$ be a locally finite graph
with no sources, and let $\tau$ be a faithful, semifinite, gauge
invariant, lower semiconitnuous trace on $C^*(E)$. Then
$(\A,\HH,\D)$ is a $QC^\infty$, $(1,\infty)$-summable, odd, local,
semifinite spectral triple (relative to $(\cn,\tilde\tau)$). For
all $a\in \A$, the operator $a(1+\D^2)^{-1/2}$ is not trace class.
If $v\in E^0$ has no sinks downstream
$$\tilde\tau_\omega(p_v(1+\D^2)^{-1/2})=2\tau(p_v),$$
where $\tilde\tau_\omega$ is any Dixmier trace associated to
$\tilde\tau$.
\end{thm}
We require the definitions of $\cn$ and $\tilde\tau$,
along with some preliminary results.

%We consider the right representation of the core $F$ on $\HH$, or
%equivalently, the (left) representation of $F^{op}$ on $\HH$.
\begin{defn} Let $End^{00}_F(X_c)$ denote the algebra of finite rank
operators on $X_c$ acting on $\HH$. Define
$\cn=(End^{00}_F(X_c))''$, and let $\cn_+$ denote the positive
cone in $\cn$.
\end{defn}

\begin{defn} Let  $T\in\cn$ and $\mu\in E^*$. Let $|v|_k=$ the
 number of paths of
 length $k$ with range $v$, and define for $|\mu|\neq 0$
$$\omega_\mu(T)=
\la S_\mu,TS_\mu\ra+\frac{1}{|r(\mu)|_{|\mu|}}\la
S_\mu^*,TS_\mu^*\ra.$$ For $|\mu|=0$, $S_\mu=p_v$, for some $v\in
E^0$, set $\omega_\mu(T)=\la S_\mu,TS_\mu\ra.$ Define
$$\tilde\tau:\cn_+\to[0,\infty],\ \ \mbox{by}\ \ \ \ \tilde\tau(T)=
\lim_{L\uparrow}\sum_{\mu\in L\subset E^*}\omega_\mu(T)$$
where $L$ is in the net of finite subsets of $E^*$.
\end{defn}

{\bf Remark} For $T,S\in\cn_+$ and $\lambda\geq 0$ we have
$$\tilde\tau(T+S)=\tilde\tau(T)+\tilde\tau(S)\ \ \
\mbox{and}\ \ \ \tilde\tau(\lambda T)=\lambda\tilde\tau(T)\ \ \mbox{where}\ \
0\times\infty=0.$$

\begin{prop}\label{tildetau}
 The function $\tilde\tau:\cn_+\to[0,\infty]$ defines a faithful normal
semifinite trace on $\cn$. Moreover,
$$End_F^{00}(X_c)\subset\cn_{\tilde\tau}:=
{\rm span}\{T\in\cn_+:\tilde\tau(T)<\infty\},$$
the domain of definition of $\tilde\tau$, and
$$\tilde\tau(\Theta^R_{x,y})=\la y,x\ra=\tau(y^*x),\ \ \ x,y\in X_c.$$
\end{prop}
\begin{proof}
First, since $\tilde\tau$ is defined as the limit of an increasing net
of sums of positive vector
functionals, $\tilde\tau$ is a positive ultra-weakly lower semicontinuous
weight on $\cn_+$, \cite{KR}, that is a normal weight. Now
 observe (using the fact that $p_v\Phi_k$ is a projection for
 all $k\in\Z$ and $v\in E^0$) that for any vertex $v\in E^0$, $k\in\Z$ and $T\in\cn_+$
\bean\tilde\tau(p_v\Phi_kTp_v\Phi_k)&=&\la \Phi_kp_v,T\Phi_kp_v\ra+
\sum_{s(\mu)=v}\la \Phi_kS_\mu,T\Phi_kS_\mu\ra\nno&+&
\sum_{r(\mu)=v}\frac{1}{|r(\mu)|_{|\mu|}}\la \Phi_kS_\mu^*,T\Phi_kS_\mu^*\ra.
\eean
If $k=0$ this is equal to $\la p_v,Tp_v\ra<\infty$. If $k>0$ we find
 \bean\tilde\tau(p_v\Phi_kTp_v\Phi_k)
 &=&\sum_{s(\mu)=v,|\mu|=k}\la S_\mu, TS_\mu\ra
 \leq\Vert T\Vert\sum_{s(\mu)=v,|\mu|=k}\tau(S_\mu^*S_\mu)\nno
 &=&\Vert T\Vert\sum_{s(\mu)=v,|\mu|=k}\tau(p_{r(\mu)})
 \leq\Vert T\Vert \tau(p_v)<\infty,\eean
 the last inequality following from the fact that $\tau$ arises from a
 graph trace, by Proposition \ref{trace=graphtrace}, and
 Equations (\ref{nosinksum}) and (\ref{sinksum}). Similarly, if $k<0$
 \hspace{-7pt}\begin{align*}\tilde\tau(p_v\Phi_kTp_v\Phi_k)
 &=\sum_{r(\mu)=v,|\mu|=|k|}
 \frac{1}{|v|_{|k|}}\la S_\mu^*, TS_\mu^*\ra
 \leq\Vert T\Vert\sum_{r(\mu)=v,|\mu|=|k|}
 \frac{1}{|v|_{|k|}}\tau(S_\mu^*S_\mu)\nno
 &=\Vert T\Vert\sum_{r(\mu)=v,|\mu|=k}
 \frac{1}{|v|_{|k|}}\tau(p_{r(\mu)})
 =\Vert T\Vert \tau(p_v)<\infty.\end{align*}

Hence $\tilde\tau$ is a finite positive function on each
$p_v\Phi_k\cn p_v\Phi_k$. Taking
limits over
finite sums of vertex projections, $p=p_{v_1}+\cdots+p_{v_n}$,
converging to the
identity, and finite sums $P=\Phi_{k_1}+\cdots+\Phi_{k_m}$,
we have for $T\in\cn_+$
$$\lim_{pP\nearrow 1}\sup\tilde\tau(pPTpP)\leq\tilde\tau(T)\leq
\lim_{pP\nearrow1}\inf\tilde\tau(pPTpP),$$ the first inequality
following from the definition of $\tilde\tau$, and the latter from
the ultra-weak lower semicontinuity of $\tilde\tau$, so for
$T\in\cn_+$
\be\lim_{pP\nearrow1}\tilde\tau(pPTpP)=\tilde\tau(T).\label{tildetaulimit}\ee
For $x\in X_c\subset\HH$, $\Theta^R_{x,x}\geq 0$ and so we compute
\bean \tilde\tau(\Theta^R_{x,x})&=& \sup_F\sum_{\mu\in F}\la
S_\mu,x(x|S_\mu)_R\ra+ \frac{1}{|r(\mu)|_{|\mu|}}\la
S_\mu^*,x(x|S_\mu^*)_R\ra\nno &=&\sup_F\sum_{\mu\in
F}\tau(\Phi(S_\mu^*x\Phi(x^*S_\mu)))+
\frac{1}{|r(\mu)|_{|\mu|}}\tau(\Phi(S_\mu
x\Phi(x^*S_\mu^*))).\eean Now since $x\in X_c$, there are only
finitely many $\omega_\mu$ which are nonzero on $\Theta^R_{x,x}$,
so this is always a finite sum, and
$\tilde\tau(\Theta^R_{x,x})<\infty$.

To compute $\Theta^R_{x,y}$, suppose that $x=S_\al S_\beta^*$ and
$y=S_\s S_\rho^*$. Then $(y|S_\mu)_R=\Phi(S_\rho S_\s^*S_\mu)$ and
this is zero unless $|\s|=|\mu|+|\rho|$. In this case, $|\s|\geq
|\mu|$ and we write $\s=\underline{\s}\overline{\s}$ where
$|\underline{\s}|=|\mu|$. Similarly, $(y|S^*_\mu)_R=\Phi(S_\rho
S_\s^*S_\mu^*)$ is zero unless $|\rho|=|\s|+|\mu|$. We also
require the computation
$$ S_\al S_\beta^* S_\rho S_\s^*S_\mu S_\mu^*=
S_\al S_\beta^* S_\rho S_\s^*\delta_{\underline{\s},\mu},\qquad
  |\s|\geq |\mu|$$
$$ S_\al S_\beta^*S_\rho S_\s^*S_\mu^*S_\mu=
S_\al S_\beta^*S_\rho S_\s^* \delta_{r(\mu),s(\s)}
\qquad|\mu|\geq|\s|.$$ Now we can compute for $|\rho|\neq|\s|$ ,
so that only one of the sums over $|\mu|=\pm(|\s|-|\rho|)$ in the
next calculation is nonempty: \bean\tilde\tau(\Theta^R_{x,y})
&=&\sum_{\mu}\tau(S_\mu^* x\Phi(y^*S_\mu))+
\sum_{\mu}\frac{1}{|r(\mu)|_{|\mu|}} \tau(S_\mu
x\Phi(y^*S_\mu^*))\nno &=&\sum_{|\mu|=|\s|-|\rho|}\tau(xy^*S_\mu
S_\mu^*)
+\sum_{|\mu|=|\rho|-|\s|}\frac{1}{|r(\mu)|_{|\mu|}}\tau(xy^*S_\mu^*
S_\mu)\nno
&=&\sum_{|\mu|=|\s|-|\rho|}\tau(xy^*\delta_{\underline{\s},\mu})
+\sum_{|\mu|=|\rho|-|\s|,r(\mu)=s(\s)}\frac{1}{|r(\mu)|_{|\mu|}}\tau(xy^*)\nno
&&\nno &=&\tau(xy^*)=\tau(y^*x)=\tau((y|x)_R)=\la y,x\ra. \eean
 When
$|\s|=|\rho|$,  we have
$$\tilde\tau(\Theta^R_{x,y})=\sum_{v\in
E^0}\tau(\Phi(p_vxy^*p_v))=\sum_{v\in E^0}\tau(y^*p_vx)$$ and the
same conclusion is obtained as above. By linearity,  whenever
$x,y\in X_c$, $\tilde\tau(\Theta^R_{x,y})=\tau((y|x)_R)$.
 For any two
$\Theta^R_{x,y}$, $\Theta^R_{w,z}\in End_F^{00}(X_c)$ we find
\bean\tilde\tau(\Theta^R_{w,z}\Theta^R_{x,y})&=&\tilde\tau(\Theta^R_{w(z|x)_R,y})
=\tau((y|w(z|x)_R)_R)
=\tau((y|w)_R(z|x)_R)\nno&=&\tau((z|x)_R(y|w)_R)
=\tilde\tau(\Theta^R_{x(y|w)_R,z})=\tilde\tau(\Theta^R_{x,y}\Theta^R_{w,z}).\eean
Hence by linearity, $\tilde\tau$ is a trace on
$End_F^{00}(X_c)\subset \cn$.

We saw previously that $\tilde\tau$ is finite on $pP\cn pP$ whenever $p$ is a
finite sum of vertex projections $p_v$ and $P$ is a finite sum of the spectral
projections $\Phi_k$.

Since $\tilde\tau$
is ultra-weakly lower semicontinuous on $pP\cn_+ pP$, it is completely additive
in the sense of \cite[Definition 7.1.1]{KR}, and therefore is normal by
\cite[Theorem 7.1.12]{KR}, which is to say, ultra-weakly continuous.

The algebra
$End^{00}_F(X_c)$ is
 strongly dense
in $\cn$, so $pPEnd^{00}_F(X_c) pP$
%\subset End_F^{00}(X_c)$
is strongly
 dense in $pP\cn pP$. Let $T\in pP\cn pP$,  and choose a bounded net
$ T_i$,  converging $*$-strongly to $T$, with $T_i\in pP End_F^{00}(X_c) pP$.
 Then, since  multiplication is jointly continuous on
bounded sets in the $*$-strong topology,
$$\tilde\tau(TT^*)=\lim_{i}\tilde\tau(T_iT_i^*)=\lim_{i}\tilde\tau(T_i^*T_i)
=\tilde\tau(T^*T).$$
Hence $\tilde\tau$ is a trace on each $pP\cn pP$ and so on $\cup_{pP} pP\cn pP$,
where the union is over all finite sums $p$ of vertex projections and finite
sums $P$ of the $\Phi_k$.

Next we want to show that $\tilde\tau$ is semifinite, so for all
$T\in\cn$ we want to find a net $R_i\geq0$ with $R_i\leq T^*T$ and
$\tilde\tau(R_i)<\infty$. Now
$$\lim_{pP\nearrow1}T^*pPT=T,\ \ \ \ T^*pPT\leq T$$
and we just need to show that $\tilde\tau(T^*pPT)<\infty$. It
suffices to show this for  $pP=p_v\Phi_k$, $v\in E^0,\ k\in\Z$. In
this case we have (with $q$ a finite sum of vertex projections and
$Q$ a finite sum of $\Phi_k$) \bean \tilde\tau(T^*p_v\Phi_kT)&=&
\lim_{qQ\nearrow 1}\tilde\tau(qQT^*p_v\Phi_kTqQ) \qquad\mbox{by
equation}\ (\ref{tildetaulimit})\nno &=&\lim_{qQ\nearrow
1}\tilde\tau(qQT^*qQp_v\Phi_kTqQ) \qquad\mbox{eventually}\
qQp_v\Phi_k=p_v\Phi_k\nno
&=&\lim_{qQ\nearrow1}\tilde\tau(qQp_v\Phi_kT^*qQTqQp_v\Phi_k)
\quad\tilde\tau\ \mbox{is a trace on}\ qQ\cn qQ\nno
&=&\lim_{qQ\nearrow1}\tilde\tau(p_v\Phi_kT^*qQTp_v\Phi_k)
=\tilde\tau(p_v\Phi_kTp_v\Phi_k)<\infty \eean Thus $\tilde\tau$ is
semifinite, normal weight on $\cn_+$, and is a trace on a dense
subalgebra. Now let $T\in\cn$. By the above
\be\tilde\tau(T^*pPT)=\tilde\tau(pPT^*TpP).\label{almost}\ee By
lower semicontinuity and the fact that $T^*pPT\leq T^*T$, the
limit of the left hand side of Equation (\ref{almost}) as $pP\to
1$ is $\tilde\tau(T^*T)$. By Equation (\ref{tildetaulimit}), the
limit of the right hand side is $\tilde\tau(TT^*)$. Hence
$\tilde\tau(T^*T)=\tilde\tau(TT^*)$ for all $T\in\cn$, and
$\tilde\tau$ is a normal, semifinite trace on $\cn$.
\end{proof}
 %and consider
%$\tilde\tau(T^*T)$. If $\tilde\tau(T^*T)$ is finite then there are  sequences
%$\{T_i\},\{S_i\}\subset End_F^{00}(X_c)_\HH$ such that $T_i\to T$ and
%$S_i\to T^*$
%(ultra weakly and so strongly) and since
%$\tilde\tau$ is ultra weakly continuous on $\tilde\tau$-finite elements and multiplication
%is jointly continuous in the strong topology
%\be\tilde\tau(T^*T)=\lim\tilde\tau(T^*_iT_i)=\lim\tilde\tau(T_iT_i^*)=
%\tilde\tau(TT^*).\label{flip}\ee
%So $\tilde\tau$ is a trace on $\tilde\tau$-finite elements. Suppose then that
%$\tilde\tau(T^*T)=\infty$. Then $\tilde\tau(TT^*)=\infty$ also, for otherwise
%we could apply the same argument as in Equation (\ref{flip}) to find that
%$\tilde\tau(T^*T)<\infty$; a contradiction. Hence $\tilde\tau:\cn\to\C$ is
%a semifinite ultra weakly lower semi-continuous trace.
%Since $\tilde\tau$ being ultra weakly lower semicontinuous is equivalent to
%normality, \cite[]{}, the proof is complete.

{\bf Notation} If $g:E^0\to\R_+$ is a faithful graph trace, we
shall write $\tau_g$ for the associated semifinite trace on
$C^*(E)$, and $\tilde\tau_g$ for the associated faithful,
semifinite, normal trace on $\cn$ constructed above.

\begin{lemma}\label{tracepvphi} Let $E$ be a locally finite graph with
no sources and a faithful graph trace $g$. Let $v\in E^0$ and
$k\in\Z$. Then
$$ \tilde\tau_g(p_v\Phi_k)\leq \tau_g(p_v)$$
with equality when $k\leq 0$ or when $k>0$ and there are no sinks
within $k$ vertices of $v$.
\end{lemma}

\begin{proof} Let $k\geq 0$. Then, by Lemma \ref{finrank} we have
\bean \tilde\tau_g\left(p_v\Phi_k\right)&=&
\tilde\tau_g\left(p_v\sum_{|\rho|=k}\Theta^R_{S_\rho,S_\rho}\right)
=\tilde\tau_g\left(\sum_{|\rho|=k}\Theta^R_{p_vS_\rho,S_\rho}\right)\nno
&=&\tau_g\left(\sum_{|\rho|=k}(S_\rho|p_vS_\rho)_R\right)=
\tau_g\left(\sum_{|\rho|=k}\Phi(S_\rho^*p_vS_\rho)\right)\nno&=&
\tau_g\left(\sum_{|\rho|=k,s(\rho)=v}S_\rho^*S_\rho\right)=
\tau_g\left(\sum_{|\rho|=k,s(\rho)=v}p_{r(\rho)}\right).\eean Now
$\tau_g(p_v)=g(v)$ where $g$ is the graph trace associated to
$\tau_g$, Proposition \ref{trace=graphtrace}, and Equation
(\ref{sinksum}) shows that\be g(v)=\sum_{|\rho|\preceq k,\
s(\rho)=v}g(r(\rho))\geq\sum_{|\rho|=k,s(\rho)=v}g(r(\rho)),\label{dblestar}\ee
with equality provided there are no sinks within $k$ vertices of
$v$ (always true for $k=0$). Hence for $k\geq 0$ we have
$\tilde\tau_g(p_v\Phi_k)\leq\tau_g(p_v),$ with equality when there
are no sinks within $k$ vertices of $v$. For $k<0$ we proceed as
above and observe that there is at least one path of length $|k|$
ending at $v$ since $E$ has no sources. Then
\begin{align}
\tilde\tau_g(p_v\Phi_{k})&=\frac{1}{|v|_k}\sum_{|\rho|=|k|,\
r(\rho)=v}\tau_g(S_\rho
p_vS_\rho^*)=\frac{1}{|v|_k}\sum_{|\rho|=|k|,\
r(\rho)=v}\tau_g(S_\rho^*S_\rho p_v)\nno
&=\frac{1}{|v|_k}\sum_{|\rho|=|k|,\
r(\rho)=v}\tau_g(p_v)=\tau_g(p_v).\qed\end{align} \hideqed
\end{proof}

\begin{prop}\label{Dixytilde=tau} Assume that the directed graph $E$ is
locally finite, has no sources and has a faithful graph trace $g$.
 For all $a\in A_c$ the operator $a(1+\D^2)^{-1/2}$ is in the
ideal $\LL^{(1,\infty)}(\cn,\tilde\tau_g)$.
\end{prop}

\begin{proof} It suffices to show that
 $a(1+\D^2)^{-1/2}\in\LL^{(1,\infty)}(\cn,\tilde\tau_g)$ for a vertex projection
$a=p_v$ for $v\in E^0$, and extending to more general $a\in A_c$
using the arguments of Lemma \ref{finrank}. Since $p_v\Phi_k$ is a
projection for all $v\in E^0$ and $k\in\Z$, we may compute the
Dixmier trace using the partial sums (over $k\in\Z$) defining the
trace of $p_v(1+\D^2)^{-1/2}$. For the partial sums with $k\geq
0$, Lemma \ref{tracepvphi} gives us \begin{equation}
\tilde\tau_g\left(p_v\sum_{0}^N(1+k^2)^{-1/2}\Phi_k\right)
\leq\sum_{k=0}^N(1+k^2)^{-1/2}\tau_g(p_v).\label{dblestarry}\end{equation}
We have equality when there are no sinks within $N$ vertices of
$v$. For the partial sums with $k<0$ Lemma \ref{tracepvphi} gives
$$\sum_{k=-N}^{-1}(1+k^2)^{-1/2}\tilde{\tau}_g(p_v\Phi_{k})=
\sum_{k=-N}^{-1}(1+k^2)^{-1/2}\tau_g(p_v),$$ and  the sequence
$$\frac{1}{\log 2N+1}\sum_{k=-N}^N(1+k^2)^{-1/2}\tilde\tau_g(p_v\Phi_k)$$
is bounded.  Hence $p_v(1+\D^2)^{-1/2}\in\LL^{(1,\infty)}$ and for
any $\omega$-limit we have
$$\tilde\tau_{g\omega}(p_v(1+\D^2)^{-1/2})=
\omega\mbox{-}\!\lim\frac{1}{\log
2N+1}\sum_{k=-N}^N(1+k^2)^{-1/2}\tilde\tau_g(p_v\Phi_k).$$ When
there are no sinks downstream from $v$, we have equality in
Equation (\ref{dblestarry}) for any $v\in E^0$ and so
$$\tilde\tau_{g\omega}(p_v(1+\D^2)^{-1/2})=2\tau_g(p_v).\qed$$
\hideqed
\end{proof}
% Computing the Dixmier trace when there are
%sinks downstream is harder.

 {\bf Remark}  Using Proposition \ref{wellbehaved}, one can
check  that
 \be res_{s=0}\tilde\tau_g(p_v(1+\D^2)^{-1/2-s})=
\frac{1}{2}\tilde\tau_{g\omega}(p_v(1+\D^2)^{-1/2}).\label{res}\ee
We will require this formula when we apply the local index
theorem.

\begin{cor}\label{compactresolvent} Assume $E$ is locally finite,
has no sources and has a faithful graph trace $g$. Then for all
$a\in A$, $a(1+\D^2)^{-1/2}\in\K_\cn$.
\end{cor}

\begin{proof} (of Theorem \ref{mainthm}.) That we
have a $QC^\infty$ spectral triple follows from Corollary
\ref{smodense}, Lemma \ref{smoalg} and Corollary
\ref{compactresolvent}. The properties of the von Neumann algebra
$\cn$ and the trace $\tilde\tau$ follow from Proposition
\ref{tildetau}. The
$(1,\infty)$-summability and the value of the Dixmier trace comes
from Proposition \ref{Dixytilde=tau}. The locality of the spectral
triple follows from Lemma \ref{smoalg}.
\end{proof}

%By Proposition \ref{wellbehaved},
%$$\lim_{s\to1/2^+}(s-1/2)\tilde\tau(p_v(1+\D^2)^{-s})=
%\frac{1}{2}\tilde\tau_\omega(p_v(1+\D^2)^{-1/2}).$$
%$$res_{s=1}\tilde\tau(p_v(1+\D^2)^{-s/2})=\tilde\tau_\omega(p_v(1+\D^2)^{-1/2})$$
%and
\section{The Index Pairing}\label{index}

Having constructed semifinite spectral triples for graph
$C^*$-algebras arising from locally finite graphs with no sources
and a faithful graph trace, we can apply the semifinite local
index theorem described in \cite{CPRS2}. See also \cite{CPRS3,CM,Hig}.

There is a $C^*$-module index, which takes its values in the
$K$-theory of the core which is described in the Appendix. The
numerical index is obtained by applying the trace $\tilde\tau$ to
the difference of projections representing the $K$-theory class.
Thus for any unitary $u$ in a matrix algebra over the graph
algebra $A$
$$\la [u],[(\A,\HH,\D)]\ra\in \tilde\tau_*(K_0(F)).$$
We compute this pairing for unitaries arising from loops (with no
exit), which provide a set of generators of $K_1(\A)$. To describe
the $K$-theory of the graphs we are considering, recall the notion
of ends introduced in Definition \ref{ends}.

\begin{lemma}\label{Kofgraph} Let $C^*(E)$ be a graph
$C^*$-algebra such that no
loop in the locally finite graph $E$ has an exit.  Then,
$$K_0(C^*(E))=\Z^{\#ends},\ \ \ \ K_1(C^*(E))=
\Z^{\#loops}.$$
\end{lemma}

\begin{proof} This follows from the continuity of $K_*$ and
\cite[Corollary 5.3]{RSz}.
\end{proof}

If $A=C^*(E)$ is nonunital, we will denote by $A^+$ the algebra
obtained by adjoining a unit to $A$; otherwise we let $A^+$ denote
$A$.

\begin{defn} Let $E$ be a locally finite graph such that $C^*(E)$ has a
faithful graph  trace $g$. Let $L$ be a loop in $E$, and denote by
$p_1,\dots,p_n$ the projections associated to the vertices of $L$
and $S_1,\dots, S_n$ the partial isometries associated to the
edges of $L$, labelled so that $S^*_nS_n=p_1$ and
$$ S^*_iS_i=p_{i+1},\ i=1,\dots,n-1,\ \ S_iS_i^*=p_i,\
i=1,\dots ,n.$$
\end{defn}
\begin{lemma}\label{loops} Let $A=C^*(E)$ be a graph $C^*$-algebra with
faithful graph trace $g$. For each loop $L$ in $E$ we obtain a
unitary in $A^+$,
$$u=1+S_{1}+S_{2}+\cdots+S_{n} -(p_1+p_2+\cdots+p_n),$$
whose $K_1$ class does not vanish. Moreover, distinct loops give
rise to distinct $K_1$ classes, and we obtain a complete set of
generators of $K_1$ in this way.
\end{lemma}
\begin{proof}The proof that $u$ is unitary is a simple computation.
 The $K_1$ class of $u$ is the generator of a copy of
$K_1(S^1)$ in $K_1(C^*(E))$, as follows from \cite{RSz}. Distinct
loops give rise to distinct copies of $K_1(S^1)$, since no loop
has an exit.
\end{proof}

\begin{prop}\label{specflow} Let $E$ be a locally finite graph with no
sources and a
faithful graph trace $g$ and $A=C^*(E)$.
 The pairing between the  spectral triple $(\A,\HH,\D)$ of
Theorem \ref{mainthm} with $K_1(A)$ is given on the generators of
Lemma \ref{loops} by
$$\la [u],[(\A,\HH,\D)]\ra=-\sum_{i=1}^n\tau_g(p_i)=-n\tau_g(p_1).$$
\end{prop}
\begin{proof}
The semifinite local index theorem, \cite{CPRS2} provides a
general formula for the Chern character of $(\A,\HH,\D)$. In our
setting it is given by a one-cochain
$$\phi_1(a_0,a_1)=res_{s=0}\sqrt{2\pi
i}\tilde\tau_g(a_0[\D,a_1](1+\D^2)^{-1/2-s}),$$ and the pairing
(spectral flow) is given by
$$sf(\D,u\D u^*)=\la [u],(\A,\HH,\D)\ra=\frac{1}{\sqrt{2\pi
i}}\phi_1(u,u^*).$$ Now $[\D,u^*]=-\sum S_{i}^*$ and
$u[\D,u^*]=-\sum_{i=1}^n p_{i}$. Using Equation (\ref{res}) and
Proposition \ref{Dixytilde=tau},
$$ sf(\D,u\D u^*)=-res_{s=0}\tilde\tau_g(\sum_{i=1}^n
p_{i}(1+\D^2)^{-1/2-s})=-\sum_{i=1}^n\tau_g(p_{i})=-n\tau_g(p_{1}),$$
the last equalities following since all the $p_{i}$ have equal
trace and there are no sinks `downstream' from any $p_i$, since no
loop has an exit.
\end{proof}

{\bf Remark} The $C^*$-algebra of the graph consisting of a single
edge and single vertex is $C(S^1)$ (we choose Lebesgue measure as
our trace, normalised so that $\tau(1)=1$). For this example, the spectral
triple we have constructed is the Dirac triple of the circle,
$(C^\infty(S^1),L^2(S^1),\frac{1}{i}\frac{d}{d\theta})$, (as can
be seen from Corollary \ref{fryingpan}.) The index theorem above
gives the correct normalisation for the index pairing on the
circle. That is, if we denote by $z$ the unitary coming from
the construction of Lemma \ref{loops} applied to this graph, then
$\la[\bar z],(\A,\HH,\D)\ra=1$.

\begin{prop}\label{C*specflow} Let $E$ be a locally finite graph with
no sources and a faithful graph trace $g$, and
 $A=C^*(E)$.  The pairing between the spectral triple
$(\A,\HH,\D)$ of Theorem \ref{mainthm} with $K_1(A)$ can be
computed as follows. Let $P$ be the positive spectral projection
for $\D$, and perform the $C^*$ index pairing of Proposition
\ref{themapH}:
$$K_1(A)\times KK^1(A,F)\to K_0(F),\ \ \ \
[u]\times[(X,P)]\to [\ker PuP]-[{\rm coker}PuP].$$ Then we have
$$sf(\D,u\D u^*)=\tilde\tau_g(\ker
PuP)-\tilde\tau_g({\rm coker}PuP)=\tilde\tau_{g*}([\ker PuP]-[{\rm
coker}PuP]).$$
\end{prop}

\begin{proof} It suffices to prove this on the generators of $K_1$ arising
from loops $L$ in $E$. Let $u=1+\sum_iS_i-\sum_ip_i$ be the
corresponding unitary in $A^+$ defined in Lemma \ref{loops}. We
will show that $\ker PuP=\{0\}$ and that
$\mbox{coker}PuP=\sum_{i=1}^n p_{i}\Phi_1$. For $a\in PX$ write
$a=\sum_{m\geq 1}a_m$. For each $m\geq 1$ write $a_m=\sum_{i=1}^n
p_{i}a_m+(1-\sum_{i=1}^n p_{i})a_m$. Then \bean
&&PuPa_m=P(1-\sum_{i=1}^n p_{i}+\sum_{i=1}^n S_{i})a_m\nno
&=&P(1-\sum^n p_{i}+\sum^n S_{i})(\sum^n p_{i}a_m)+ P(1-\sum^n
p_{i}+\sum^n S_{i})(1-\sum^n p_{i})a_m\nno &=&P\sum^n S_{i}a_m
+P(1-\sum^n p_{i})a_m\nno &=&\sum^n S_{i}a_m+(1-\sum^n
p_{i})a_m.\eean It is clear from this computation that $PuPa_m\neq
0$ for $a_m\neq 0$.

Now suppose $m\geq 2$. If $\sum_{i=1}^n p_{i}a_m=a_m$ then
$a_m=\lim_N\sum^N_{k=1} S_{\mu_k}S_{\nu_k}^*$ with
$|\mu_k|-|\nu_k|=m\geq 2$ and $S_{{\mu_k}_1}=S_{i}$ for some $i$.
So we can construct $b_{m-1}$ from $a_m$ by removing the initial
$S_{i}$'s. Then $a_m=\sum_{i=1}^n S_{i}b_{m-1}$, and
$\sum_{i=1}^np_{i}b_{m-1}=b_{m-1}$. For arbitrary $a_m$, $m\geq
2$, we can write $a_m=\sum_ip_ia_m+(1-\sum_ip_i)a_m$, and so \bean
a_m&=&\sum^n p_{i}a_m+(1-\sum^n p_{i})a_m\nno &=&\sum^n
S_{i}b_{m-1}+(1-\sum^n p_{i})a_m\ \ \ \mbox{and by adding
zero}\nno &=&\sum^n S_{i}b_{m-1}+(1-\sum^n
p_{i})b_{m-1}+\bigl(\sum^n S_{i}+(1-\sum^n p_{i})\bigr)(1-\sum^n
p_{i})a_m\nno
&=&ub_{m-1}+u(1-\sum^np_i)a_m\nno&=&PuPb_{m-1}+PuP(1-\sum^n
p_{i})a_m.\eean Thus $PuP$ maps onto $\sum_{m\geq2}\Phi_mX$.

For $m=1$, if we try to construct $b_0$ from $\sum_{i=1}^n
p_{i}a_1$ as above, we find $PuPb_0=0$ since $Pb_0=0$. Thus
$\mbox{coker}PuP=\sum^n p_{i}\Phi_1X$. By Proposition
\ref{specflow}, the pairing is then
\begin{align} sf(\D,u\D u^*)&=-\sum^n\tau_g(p_{i})=-\tilde\tau_g(\sum^n
p_{i}\Phi_1)\nno
&=-\tilde\tau_{g*}([\mbox{coker}PuP])=-\tilde\tau_g(\mbox{coker}PuP).\end{align}
Thus we can recover the numerical index using $\tilde\tau_g$ and
the $C^*$-index.
\end{proof}

The following example shows that the semifinite index provides
finer invariants of directed graphs than those obtained from the
ordinary index. The ordinary index computes the
pairing between the $K$-theory and $K$-homology of  $C^*(E)$,
while the semifinite index also depends on the core and the gauge action.

\begin{cor}[Example]\label{fryingpan} Let $C^*(E_n)$ be the algebra
determined by the graph 
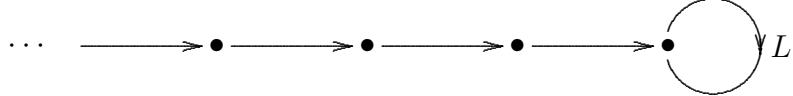
\vspace{-2pt}
\[
\beginpicture

\setcoordinatesystem units <1cm,1cm>

\setplotarea x from 0 to 12, y from -0.5 to 0.5

\put{$\cdots$} at 0.5 0

\put{$\bullet$} at 3 0

\put{$\bullet$} at 5 0

\put{$\bullet$} at 7 0

\put{$\bullet$} at 9 0

\put{$L$} at 10.5 0

\circulararc -325 degrees from 9 0.2 center at 9.6 0

\arrow <0.25cm> [0.2,0.5] from 1.2 0 to 2.8 0

\arrow <0.25cm> [0.2,0.5] from 3.2 0 to 4.8 0

\arrow <0.25cm> [0.2,0.5] from 5.2 0 to 6.8 0

\arrow <0.25cm> [0.2,0.5] from 7.2 0 to 8.8 0

\arrow <0.25cm> [0.2,0.5] from 10.228 0.1 to 10.226 -0.1

\endpicture
\]
\smallskip
where the loop $L$ has $n$ edges. Then $C^*(E_n)\cong
C(S^1)\otimes\K$ for all $n$, but $n$ is an invariant of the pair
of algebras $(C^*(E_n),F_n)$ where $F_n$ is the core of
$C^*(E_n)$.
\end{cor}
\begin{proof} Observe that the
graph $E_n$ has a one parameter family of faithful graph traces,
specified by $g(v)=r\in \R_+$ for all $v\in E^0$.

First consider the case where the graph consists only of the loop
$L$.  The $C^*$-algebra $A$ of this graph is isomorphic to
$M_n(C(S^1))$, via
$$ S_i\to e_{i,i+1},\ i=1,\dots,n-1,\ \ S_n\to id_{S^1}e_{n,1},$$
where the $e_{i,j}$ are the standard matrix units for $M_n(\C)$,
\cite{aH}. The unitary
$$S_1S_2\cdots S_n+S_2S_3\cdots S_1+\cdots+S_nS_1\cdots S_{n-1}$$ is
mapped to the orthogonal sum $id_{S^1}e_{1,1}\oplus
id_{S^1}e_{2,2}\oplus\cdots\oplus id_{S^1}e_{n,n}$. The core $F$
of $A$ is $\C^n=\C[p_1,\dots,p_n]$. Since $KK^1(A,F)$ is equal to
$$\oplus^n
KK^1(A,\C)=\oplus^nKK^1(M_n(C(S^1)),\C)=\oplus^nK^1(C(S^1))=\Z^n$$
we  see that $n$ is the rank of $KK^1(A,F)$ and so an invariant,
but let us link this to the index computed in Propositions
\ref{specflow} and \ref{C*specflow} more explicitly. Let
$\phi:C(S^1)\to A$ be given by $\phi(id_{S^1})=S_1S_2\cdots
S_n\oplus \sum_{i=2}^ne_{i,i}$. We observe that
$\D=\sum_{i=1}^np_i\D$ because the `off-diagonal' terms are $p_i\D
p_j=\D p_ip_j=0$. Since $S_1S_1^*=S^*_nS_n=p_1$, we find (with $P$
the positive spectral projection of $\D$)
$$\phi^*(X,P)=(p_1X,p_1Pp_1)\oplus\mbox{degenerate\ module}\in KK^1(C(S^1),F).$$
Now let $\psi:F\to\C^n$ be given by
$\psi(\sum_jz_jp_j)=(z_1,z_2,...,z_n)$. Then
$$\psi_*\phi^*(X,P)=\oplus_{j=1}^n(p_1Xp_j,p_1Pp_1)\in\oplus^nK^1(C(S^1)).$$
Now $X\cong M_n(C(S^1))$, so $p_1Xp_j\cong C(S^1)$ for each
$j=1,\dots,n$. It is easy to check that $p_1\D p_1$ acts by
$\frac{1}{i}\frac{d}{d\theta}$ on each $p_1Xp_j$, and so our
Kasparov module maps to
$$\psi_*\phi^*(X,P)=\oplus^n(C(S^1),P_{\frac{1}{i}\frac{d}{d\theta}})\in
\oplus^nK^1(C(S^1)),$$ where $P_{\frac{1}{i}\frac{d}{d\theta}}$ is
the positive spectral projection of
$\frac{1}{i}\frac{d}{d\theta}$. The pairing with $id_{S^1}$ is
nontrivial on each summand, since $\phi(id_{S^1})=S_1\cdots
S_n\oplus \sum_{i=2}^ne_{i,i}$ is a unitary mapping $p_1Xp_j$ to
itself for each $j$. So  we have, \cite{HR},
\begin{align}id_{S^1}\times\psi_*\phi^*(X,P)&=\sum^n_{j=1}Index(Pid_{S^1}P:p_1PXp_j\to
p_1PXp_j)\nno &=-\sum_{j=1}^n[p_j]\in K_0(\C^n).\end{align} By
Proposition \ref{C*specflow}, applying the trace to this index
gives $-n\tau_g(p_1)$. Of course in Proposition \ref{C*specflow}
we used the unitary $S_1+S_2+\cdots+S_n$, however in $K_1(A)$
$$[S_1S_2\cdots S_n]=[S_1+S_2+\cdots+S_n]=[id_{S^1}].$$
To see this, observe that
$$(S_1+\cdots+S_n)^n=S_1S_2\cdots S_n+S_2S_3\cdots
S_1+\cdots+S_nS_1\cdots S_{n-1}.$$ This is the orthogonal sum of
$n$ copies of $id_{S^1}$, which is equivalent in $K_1$ to
$n[id_{S^1}]$. Finally, $[S_1+\cdots+S_n]=[id_{S^1}]$ and so
$$[(S_1+\cdots+S_n)^n]=n[S_1+\cdots+S_n]=n[id_{S^1}].$$
Since we have cancellation in $K_1$, this implies that the class
of $S_1+\cdots+S_n$ coincides with the class of $S_1S_2\cdots
S_n$.

Having seen what is involved, we now add the infinite path on the
left. The core becomes $\K\oplus\K\oplus\cdots\oplus\K$ ($n$
copies). Since $A=C(S^1)\otimes\K= M_n(C(S^1))\otimes\K$, the
intrepid reader can go through the details of an argument like the
one above, with entirely analogous results.
\end{proof}

Since the invariants obtained from the semifinite index are finer
than the isomorphism class of $C^*(E)$, depending as they do on
 $C^*(E)$ and the gauge action, they can be regarded as invariants of the
differential structure. That is, the core $F$ can be recovered
from the gauge action, and we regard these invariants as arising
from the differential structure defined by $\D$. Thus in this
case, the semifinite index produces invariants of the differential
topology of the noncommutative space $C^*(E)$.

 \appendix
\vspace{-10pt}
\section{Toeplitz Operators on $C^*$-modules}
\vspace{-7pt} In this Appendix we define a bilinear product
$$ K_1(A)\times KK^1(A,B)\to K_0(B).$$
Here we suppose that $A, B$ are ungraded $C^*$-algebras. This
product should be the Kasparov product,
 though it is difficult to compare the two (see the footnote to
Proposition \ref{themapH} below).

We denote by $A^+$ the minimal (one-point) unitization if $A$ is
nonunital. Otherwise $A^+$ will mean $A$. To deal with unitaries
in matrix algebras over $A$, we recall that $K_1(A)$ may be
defined by considering unitaries in matrix algebras over $A^+$
which are equal to $1_n$ mod $A$ (for some $n$), \cite[p 107]{HR}.

We consider odd Kasparov $A$-$B$-modules. So let $E$ be a fixed
countably generated ungraded $B$-$C^*$-module, with $\phi:A\to
End_B(E)$ a $*$-homomorphism, and let $P\in End_B(E)$ be such that
$a(P-P^*), a(P^2-P), [P,a]$ are all compact endomorphisms. Then by
\cite[Lemma 2, Section 7]{K}, the pair $(\phi,P)$ determines a
$KK^1(A,B)$ class, and every class has such a representative. The
equivalence relations on pairs $(\phi,P)$ that give $KK^1$ classes
are unitary equivalence $(\phi,P)\sim (U\phi U^*,UPU^*)$ and
homology, $P_1\sim P_2$ if $P_1\phi_1(a)-P_2\phi_2(a)$ is a
compact endomorphism for all $a\in A$.

Now let $u\in M_m(A^+)$ be a unitary, and $(\phi,P)$ a
representative of a $KK^1(A,B)$ class. Observe that $(P\otimes
1_m)E\otimes\C^m$ is a $B$-module, and so can be extended to a
$B^+$ module. Writing $P_m=P\otimes 1_m$, the operator
$P_m\phi(u)P_m$ is Fredholm, since (dropping the $\phi$ for now)
$$ P_muP_m P_mu^*P_m=P_m[u,P_m]u^*P_m+P_m,$$
and this is $P_m$ modulo compact endomorphisms. To ensure that
$\ker P_muP_m$ and $\ker P_mu^*P_m$ are closed submodules, we need
to know that $P_muP_m$ is regular, but by \cite[Lemma 4.10]{GVF},
we can always replace $P_muP_m$ by a regular operator on a larger
module. Then the index of $P_muP_m$ is defined as the index of
this regular operator, so there is no loss of generality in
supposing that $P_muP_m$ is regular. Then we can define
$$Index(P_muP_m)=[\ker P_muP_m]-[\mbox{coker} P_muP_m]\in K_0(B).$$
This index lies in $K_0(B)$ rather than $K_0(B^+)$ by
\cite[Proposition 4.11]{GVF}.
So given $u$ and $(\phi,P)$ we  define a $K_0(B)$ class  by setting
$$ u\times (\phi,P)\to [\ker P_muP_m]-[\mbox{coker} P_muP_m].$$
Observe the following. If $u=1_m$ then
$1_m\times (\phi,P)\to Index(P_m)=0$
so for any $(\phi,P)$ the map defined on unitaries sends the
identity to zero. Given the unitary $u\oplus v\in M_{2m}(A^+)$
(say) then
$$ u\oplus v\times (\phi,P)\to Index(P_{2m}(u\oplus
v)P_{2m})=Index(P_muP_m)+Index(P_mvP_m),$$ so for each $(\phi,P)$
the map respects direct sums. Finally, if $u$ is homotopic through
unitaries to $v$, then $P_muP_m$ is norm  homotopic
 to $P_mvP_m$, so
$$ Index(P_muP_m)=Index(P_mvP_m).$$
 By the universal property of $K_1$,
\cite[Proposition 8.1.5]{RLL}, for each $(\phi,P)$ as above there
exists a unique homomorphism $H_P:K_1(A)\to K_0(B)$ such that
$$ H_P([u])=Index(P_muP_m).$$
Now observe that $H_{UPU^*,U\phi(\cdot)U^*}=H_{P,\phi}$ since
$$Index(UPU^*(U\phi(u)U^*)UPU^*)=Index(UPuPU^*)=Index(PuP).$$
 The homomorphisms $H_P$ are bilinear, since \bean
H_{P\oplus Q}([u])&=&Index((P\oplus
Q)(\phi(u)\oplus\psi(u))(P\oplus Q))\nno
&=&Index(P\phi(u)P)+Index(Q\psi(u)Q)=H_P([u])+H_Q([u]).\eean
Finally, if $(\phi_1,P_1)$ and $(\phi_2,P_2)$ are homological, the
classes defined by $(\phi_1\oplus\phi_2,P_1\oplus 0)$ and
$(\phi_1\oplus\phi_2,0\oplus P_2)$ are operator homotopic, \cite[p
562]{K}, so \bean Index(P_1\phi_1(u)P_1)&=&Index((P_1\oplus
0)(\phi_1(u)\oplus\phi_2(u))(P_1\oplus 0))\nno &=&Index((0\oplus
P_2)(\phi_1(u)\oplus\phi_2(u))(0\oplus P_2))\nno
&=&Index(P_2\phi_2(u)P_2).\eean So  $H_P$ depends
only on the $KK$-equivalence class of $(\phi,P)$. Thus
\begin{prop}\label{themapH} With the notation above, the map\footnote{As noted
at the end of the introduction, Nigel Higson has shown us a proof that the
map $H$ is equal to the Kasparov product. The Kasparov module defined by
$PuP$ in $KK^0(\C,B)=K_0(B)$ is not a product Kasparov module, but the
class of the product of representatives $u,P$ coincides with the
class of $PuP$. }
$$H:K_1(A)\times KK^1(A,B)\to K_0(B)$$
$$H([u],[(\phi,P)]):=[\ker(PuP)]-[{\rm coker}PuP]$$
is bilinear.
\end{prop}
This is a kind of spectral flow, where we are counting the net
number of eigen-$B$-modules which cross zero along any path from
$P$ to $uPu^*$.

% \section{}
% \label{}


\vspace{-10pt}
\begin{thebibliography}{00000}

\bibitem{aH} A. an Huef, Honours thesis, University of
Newcastle, 1994

\bibitem{BPRS} T. Bates, D. Pask, I. Raeburn, W. Szymanski,
{\em The $C^*$-Algebras of Row-Finite Graphs}, New York J. Maths
{\bf 6} (2000) pp 307-324

\bibitem{BR} O. Bratteli, D. Robinson, {\em Operator Algebras and
Quantum Statistical Mechanics 1}, Springer-Verlag, 2nd Ed, 1987

\bibitem{CPS2} A. Carey, J. Phillips, F. Sukochev,
{\em Spectral Flow and Dixmier Traces},
 Advances in Mathematics, {\bf 173} (2003) pp 68-113

%\bibitem{CPRS} A. Carey, J. Phillips, A. Rennie, F. Sukochev,
%{\em The Hochschild Class of the Chern Character for
%Semifinite Spectral Triples},
%su

\bibitem{CPRS1} A. Carey, J. Phillips, A. Rennie, F. Sukochev, {\em The Hochschild Class of the Chern Character of
Semifinite Spectral Triples}, Journal of Functional Analysis, {\bf
213} (2004) pp 111-153


\bibitem{CPRS2} A. Carey, J. Phillips, A. Rennie, F. Sukochev,
{\em The Local Index Theorem in Semifinite von Neumann Algebras I:
Spectral Flow}, to appear in Advances in Mathematics

\bibitem{CPRS3} A. Carey, J. Phillips, A. Rennie, F. Sukochev,
{\em The Local Index Theorem in Semifinite von Neumann Algebras
II: The Even Case}, to appear in Advances in Mathematics

 \bibitem{C} A. Connes,
        {\em Noncommutative Geometry},
        Academic Press, 1994

        %\bibitem[C1]{C1} A. Connes, {\em Gravity Coupled with
        %Matter and the Foundation of Noncommutative Geometry},
        %Commun. Math. Phys. {\bf 182} (1996), 155--176.

       % \bibitem[C2]{C2} A. Connes, {\em Cyclic Cohomology,
       % Noncommutative Geometry and Quantum Group Symmetries},
        %Lecture Notes in Math., 1831, Springer, Berlin, (2004) pp
        %1-71

        \bibitem{CM} A. Connes, H. Moscovici,
{\em The Local Index Formula in Noncommutative Geometry}, GAFA {\bf
5} (1995) 174-243

\bibitem{DHS} K. Deike, J. H. Hong, W. Szymanski, {\em Stable
Rank of Graph Algebras. Type I Graph Algebras and Their Limits},
Indiana. Univ. Math. J. {\bf 52} No. 4 (2003) pp 963-979

\bibitem{Dix} J. Dixmier, {\em Von Neumann Algebras},
North-Holland, 1981


\bibitem{FK} T. Fack and H. Kosaki, \emph{Generalised $s$-numbers of
$\tau$-measurable operators}, Pacific J. Math. {\bf 123} (1986),
269--300

\bibitem{GGISV} V. Gayral, J.M. Gracia-Bond\'{i}a,
B. Iochum, T. Sch\"{u}cker, J.C. Varilly, {\em Moyal Planes are
Spectral Triples}, Comm. Math. Phys. {\bf 246} (2004) pp 569-623

\bibitem{GVF} J. M. Gracia-Bond\'{i}a, J. C. Varilly, H. Figueroa,
{\em Elements of Non-commutative Geometry}, Birkhauser, Boston, 2001

\bibitem{Hig} N. Higson, {\em The Local Index Formula in
Noncommutative Geometry}, Contemporary Developments in Algebraic
$K$-Theory, ictp Lecture Notes, no 15, (2003), pp 444-536

\bibitem{HR} N. Higson, J. Roe, {\em Analytic $K$-Homology},
Oxford University Press, 2000

\bibitem{H} J. v.B. Hjelmborg, {\em Purely Infinite and Stable
$C^*$-Algebras of Graphs and Dynamical Systems}, Ergod. Th. \&
Dynam. Sys. {\bf 21} (2001), pp 1789-1808


\bibitem{KR} R.V. Kadison, J. R. Ringrose, {\em Fundamentals of the Theory of
Operator Algebras. Vol II Advanced Theory}, Academic Press, 1986


\bibitem{K} G. G. Kasparov, {\em The Operator $K$-Functor and
Extensions of $C^*$-Algebras}, Math. USSR. Izv. {\bf 16} No. 3
(1981), pp 513-572

\bibitem{kpr} A. Kumjian, D. Pask and I. Raeburn, {\em Cuntz-Krieger
algebras of directed graphs}, Pacific J. Math. {\bf 184} (1998),
161--174.

\bibitem{KPRR} A. Kumjian, D. Pask, I. Raeburn, J. Renault,
{\em Graphs, Groupoids and Cuntz-Krieger Algebras}, J. Funct.
Anal. {\bf 144} (1997) pp 505-541

\bibitem{L} E. C. Lance,
{\em Hilbert $C^*$-Modules}, Cambridge University Press, Cambridge,
1995


\bibitem{mal} A. Mallios,
{\em Topological Algebras, Selected Topics}, Elsevier Science
Publishers B.V., 1986

%\bibitem[P]{P} R. J. Plymen, {\em Strong Morita equivalence,
%spinors and symplectic spinors},
%J. Operator Theory {\bf 16} (1986), pp 305-324.

\bibitem{PR} D. Pask, I. Raeburn, {\em On the K-Theory of
Cuntz-Krieger Algebras}, Publ. RIMS, Kyoto Univ., {\bf 32} No. 3
(1996) pp 415-443

\bibitem{PhR} J. Phillips, I. Raeburn, {\em An Index Theorem for Toeplitz
Operators with Noncommutative Symbol Space}, J. Funct. Anal. {\bf 120} no. 2
 (1994)  pp 239-263

\bibitem{PRen} D. Pask, A. Rennie, {\em One Dimensional Noncommutative
Manifolds from  Graph $C^*$-Algebras}, in
preparation

\bibitem{RW} I. Raeburn and D. P. Williams, {\em Morita Equivalence and
Continuous-Trace $C^*$-Algebras}, Math. Surveys \& Monographs,
vol. 60, Amer. Math. Soc., Providence, 1998.

%\bibitem[RS]{RS} M. Reed and B. Simon,
%{\em Volume I: Functional Analysis, Volume II: Fourier Analysis,
%Self-Adjointness}, Academic Press, 1980

\bibitem{RSz} I. Raeburn, W. Szymanski, {\em Cuntz-Krieger
Algebras of Infinite Graphs and Matrices}, Trans. Amer. Math. Soc.
{\bf 356} no. 1 (2004) pp 39-59

\bibitem{R} I. Raeburn, {\em Graph Algebras: $C^*$-Algebras we can
see}, CBMS Lecture Notes, to appear

\bibitem{R1} A. Rennie,
{\em Smoothness and Locality for Nonunital Spectral Triples},
$K$-theory, {\bf 28}(2) (2003) pp 127-165

\bibitem{R2} A. Rennie,
{\em Summability for Nonunital Spectral Triples}, $K$-theory, {\bf
31} (2004) pp 71-100

%\bibitem[RV]{RV} A. Rennie, J. Varilly, {\em Connes' Spin Manifold
%Theorem}, preprint

\bibitem{RLL} M. R\o rdam, F. Larsen, N. J. Laustsen, {\em An
Introduction to $K$-Theory and $C^*$-Algebras}, LMS Student Texts,
49, CUP, 2000

\bibitem{LBS} Larry B. Schweitzer,
{\em A Short Proof that $M_n(A)$ is local if $A$ is Local and
Fr\'{e}chet}, Int. J. math. {\bf 3} No.4 581-589 (1992)

\bibitem{SZ} S. Str\u{a}til\u{a}, L. Zsid\'{o}, {\em
Lectures on von Neumann Algebras}, Abacus Press, 1975

\bibitem{T} Mark Tomforde, {\em Real Rank Zero and Tracial States of
$C^*$-Algebras Associated to Graphs}, math.OA/0204095 v2

% notes:

% \bibitem{label2}
% If there is a note, it should come last:
% \bibitem{label3} \note
% \end{subbibitems}

%\bibitem{}

\end{thebibliography}
\end{document}